\documentclass[12pt,a4paper]{article}
\usepackage[T1]{fontenc}
\usepackage[utf8]{inputenc}
\usepackage{lmodern}
\usepackage{textcomp}
\usepackage[english]{babel}
\usepackage{booktabs}
\usepackage{graphicx}
\usepackage{abstract}
\usepackage{indentfirst, comment}
\usepackage[colorlinks=true,linkcolor=blue,citecolor=blue,urlcolor=blue]{hyperref}
\usepackage{authblk}
\usepackage{appendix}
\usepackage{algorithmic}
\usepackage{lscape}
\usepackage[ruled, vlined]{algorithm2e}
\usepackage{multirow}
\usepackage{wrapfig}
\usepackage{subfigure}
\usepackage{float}
\usepackage{subfloat}
\usepackage{natbib}
\usepackage{fullpage}
 \usepackage{pgfplots}
\usepackage{xcolor}
\usepackage{xfrac}
\usepackage{rotating}
\usepackage{comment}
\SetKwRepeat{Do}{do}{while}%
\usepackage{amsthm,amsfonts,amsmath,amstext,amssymb,amsopn}
\usepackage{comment}
\newcommand{\Depot}{\includegraphics[scale=1]{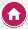}}%
\newcommand{\Pickup}{\includegraphics[scale=1]{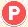}}%
\newcommand{\Delivery}{\includegraphics[scale=1]{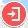}}%

\usepackage{tikz}
\usepackage{setspace}
\usepackage{pgf}

\usetikzlibrary{spy}
\usetikzlibrary{arrows}
\usetikzlibrary{patterns,decorations.pathreplacing}
\usetikzlibrary{intersections}
\usepackage[font=small,labelfont=bf,labelsep=period,tableposition=top]{caption}
\usepackage{natbib}
\pgfplotsset{compat=1.12}
\usepackage{xcolor}
\usepackage{wrapfig}
\usepackage{framed}
\usepackage{pgfplots}
\usepackage{graphicx}
\usepackage{eurosym}
\usepackage{pgf-pie}
\usepackage{soul}
\usepackage{ulem}
\usepackage{float}
\usepackage{subfloat}
\usepackage{comment}

\definecolor{dark-gray}{gray}{0.1}
\newcommand{\TCDARP}{T-CDARP}
\newcommand{\CCDARP}{C-CDARP}
\newcommand{\TCCDARP}{TC-CDARP}
\newcommand{\CDARP}{UC}
\newcommand{\NCDARP}{NC}
\newcommand{\GAP}{\textit{GAP}}
\newcommand{\SAV}{\textit{SAV}}

\title{Optimization models for fair horizontal collaboration in demand-responsive transportation}
\author[1]{E. Angelelli}
\author[2]{V. Morandi}
\author[1]{M.G. Speranza}

\affil[1]{University of Brescia}
\affil[2]{Free University of Bolzano/Bozen}

\begin{document}

\maketitle

	\begin{abstract}
	
The  advances in information and communication technology are changing the way people move. Companies that offer demand-responsive transportation services have the opportunity to reduce their costs and increase their revenues through collaboration, while at the same time reducing the environmental impact of their operations. We consider the case of companies, offering a shared taxi service, that are involved in horizontal collaboration and present mixed integer programming models for the optimization of their routes that embed constraints aimed at balancing the workload exchange.
 These constraints bound the imbalance in terms of traveled time  and/or served customers  to be less than thresholds agreed in advance by the companies. We also present a heuristic algorithm and show the benefits of the collaboration.

		\end{abstract}

\section{Introduction}
			
Conventional transportation systems typically compromise between coverage of the demand and cost of the service. The service is usually based on fixed routes and fixed schedules.
Good service can be provided in areas characterized by high-density demand but  cannot be afforded in others, where low-frequency service may be cost-ineffective and leave  users unsatisfied.
Taxis, on the other hand, offer a good but expensive service.

The advances in information and communication technologies have created  opportunities for new transportation services.
Demand-Responsive Transport (DRT) systems have been suggested as systems capable to satisfy the dynamic nature of the demand (see, for example, \cite{amirgholy2016demand}).
The term DRT has been increasingly used to refer to systems aimed at feeding, by means of small buses or vans or taxis, conventional transportation systems.
Over time, the concept of DRT  has been broadened to include a  range of flexible, demand-responsive transport services, that neither rely on fixed routes nor on fixed schedules.
Another term that has become very popular, and often used as a synonym of DRT systems, is that of Flexible Transport Services (FTSs) (see, for example, \cite{mulley2009flexible}).
In \cite{atasoy2015concept} a  system, called flexible mobility on demand, was introduced which provides different levels of service to each customer request.
The system makes use of three services, taxi, shared taxi, and minibus. In the shared taxi service customers accept small detours from their original direct path and share part of their ride with others, while the minibus service works as a regular bus service with fixed schedules. In \cite{martinez2017assessing} the impact of the introduction of  a shared taxi system and of a dynamic bus-like service with minibuses, where customers pre-book their service at least 30 minutes in advance and walk short distances to a designated stop, is explored.

Optimization of a shared taxi system is closely related to the stream of literature on Dial-A-Ride Problems (DARPs). The class of DARPs covers optimization models aimed at finding the best possible routes for vehicles that pick up from home and deliver to destination people that can share a portion of their trip with other passengers. 
The DARPs are
pick-up and delivery vehicle routing problems with time windows and additional constraints as on the maximum duration of each passenger trip.
The original DARP formulation dates back to \cite{cordeau2006branch} and the most recent heuristic is provided in \cite{gschwind2019adaptive}.
%
%
Several specialized DARPs have been studied, as non-profit DARPs for elderly and disabled people, airport shuttle services, healthcare services and public transportation. We refer the reader to  \cite{molenbruch2017typology} for a recent survey.

%
%


DRT services in big cities may be offered by different competing companies, either public or private.
%
%
Such companies, like many others, may benefit from horizontal collaboration which has a potential to reduce costs and increase profits.
Horizontal collaboration has been exploited in the last decades by several companies, but not as many as one might have expected.
Collaboration allows transportation companies to better exploit the vehicle capacity, with the desirable side effect of reducing the CO2 emissions.
However, creating a collaboration scheme that is beneficial to all the companies involved is far from being straightforward.
The results of a large-scale survey on the potential benefits of and impediments for horizontal collaboration among logistic companies in Flanders were presented in \cite{cruijssen2007horizontal}.
In \cite{basso2019survey} a wide discussion of the different forms of collaboration in logistics and of the obstacles and difficulties of practical implementation of horizontal collaboration was provided while in \cite{ALOUI2021100291} a more general systematic literature review on collaborative sustainable transportation was offered.
%
%
As discussed in \cite{gansterer2018collaborative} and in \cite{gansterer2020shared}, finding  mechanisms that make horizontal collaboration beneficial to the involved companies is a topic of great interest for the transportation community.


In the literature, approaches based on game theory were also proposed.
For example, a linear model was used in \cite{lozano2013cooperative} to study the cost savings that different companies may achieve when they merge their transportation requirements.
In \cite{ozener2008allocating} a logistics network was studied in which shippers collaborate and bundle their shipment requests to negotiate better rates with a common carrier and cost-allocation mechanisms were proposed to ensure the sustainability of the collaboration.
%
%
In \cite{gansterer2021prisoners} the prisoners' dilemma was applied in a collaborative pick-up and delivery setting in which requests are traded using combinatorial auction mechanisms.


In the optimization field, only a few papers have proposed  models aimed at optimizing the routes of the companies involved in a collaboration initiative.
An  optimization model for an arc routing problem was analyzed in \cite{fernandez2016collaboration} to model collaboration in truckload shipping.
A lower bound on the individual profit of each carrier is set in the optimization model to guarantee that all carriers benefit from the collaboration.
A collaborative version of a routing problem with profits was  proposed in \cite{defryn2016selective}, where the customers of different companies can be served with the joint available resources of companies joining the coalition.
The formulations incorporate different cost allocation rules that model the desired behavior of the participants.
In \cite{fernandez2018shared} a routing problem is optimized where the companies that collaborate share only some of their customers.
In \cite{gansterer2020assignment} a pick-up and delivery collaborative problem approaching fairness among companies is provided.
They model the problem so that each company is guaranteed to serve at least a minimum amount of its own customers, at least a fixed portion of other companies' customers while profit is bounded to be greater or equal than the profit obtained in a non-collaborative setting.

%
The collaboration principles have been applied to dial-a-ride service providers in \cite{molenbruch2017benefits}, where companies share their customers as if they were a single company trying to exploit the overall available resources in order to optimize the overall profit. Afterward, the profit is distributed on the basis of the number of customers each company has served. The drawback of this approach is that some companies, especially small companies, may serve only few customers which will make the collaboration not sustainable.
 As stated in \cite{cruijssen2007horizontal}, service providers strongly believe in the potential benefits of horizontal collaboration. However, the authors reported that a great majority of companies believe that smaller companies may lose customers or get pushed out of the market. Thus, the authors recommended to explicitly take these issues into account.

In this paper we consider the case of  shared taxi service providers involved in horizontal collaboration by taking explicitly into account the issue of fairness towards each company.
We propose a collaboration scheme where each  company is allowed to serve customers of other companies, if beneficial. The exchange  is, however, constrained to make the collaboration scheme acceptable to all companies, in particular to those of small or relatively small size.
The  optimization models we propose optimize the overall system cost while embedding constraints aimed at balancing the exchanged workload in terms of traveling time and/or number of customers served.
In order to balance the exchanged time workload, for each company we consider the time spent in serving other companies' customers and the time spent by other companies in serving its customers; the difference between these two quantities is  bounded to be as small as desired so that the workload remains similar to the one the company would experience in a non-collaborative setting.
Observe that it does not matter how many customers a company concedes to or acquires from others, as the bound is on the absolute value of the time balance, i.e. the  exchanged time workload, for each individual company.
%
%
As companies aim at maintaining or increasing their customer base, in a collaboration initiative they may be interested to serve  as many customers as they would have served in a non-collaborative setting. Thus, we also consider an optimization model where the imbalance of each company, in number of served customers, is  bounded. An optimization model where both the workload time exchange and the customers exchange are bounded is also studied.


%
An Adaptive Large Neighbourhood Search (ALNS) algorithm is presented for the solution of the proposed optimization models, for which we implemented a number of destroy/repair operators taken from the literature and two new destroy and one new repair operators which proved to be effective in finding good quality solutions.
%
%
Extensive computational tests are performed on a set of 112 instances randomly generated for this problem setting in order to assess the performance of the ALNS heuristic and the effectiveness of the proposed models. The results show that the average error generated by the algorithm, when compared with the optimal solution, is about  0.3\%. 
It is also shown that the optimization models  guarantee to each company a high level of fairness while attaining a total routing cost which is very close to the one obtained for the case no bound is set on the time and customer exchange.
Savings up to 30\% are obtained through fair collaboration, only slightly less than the savings that might be obtained through unconstrained, and thus possibly highly unbalanced, collaboration.
Finally, we will show that, in a multi-period horizon scenario, the imbalance of the companies on a single day may be used to offset the imbalance of the day after in order to control the long run balance.

The remainder of the paper is organized as follows. In Section \ref{mot}, we illustrate a motivating example. In Section \ref{models}, we introduce and discuss the collaborative formulations proposed. In Section \ref{heur}, we present the ALNS heuristic. In Section \ref{result}, we discuss the results of the computational study. Finally, in Section \ref{conclusions}, we present some concluding remarks.

\subsection{Motivating example}\label{mot}

\begin{figure}
\centering
\subfigure[Before routing ]
{\includegraphics[width=7cm]{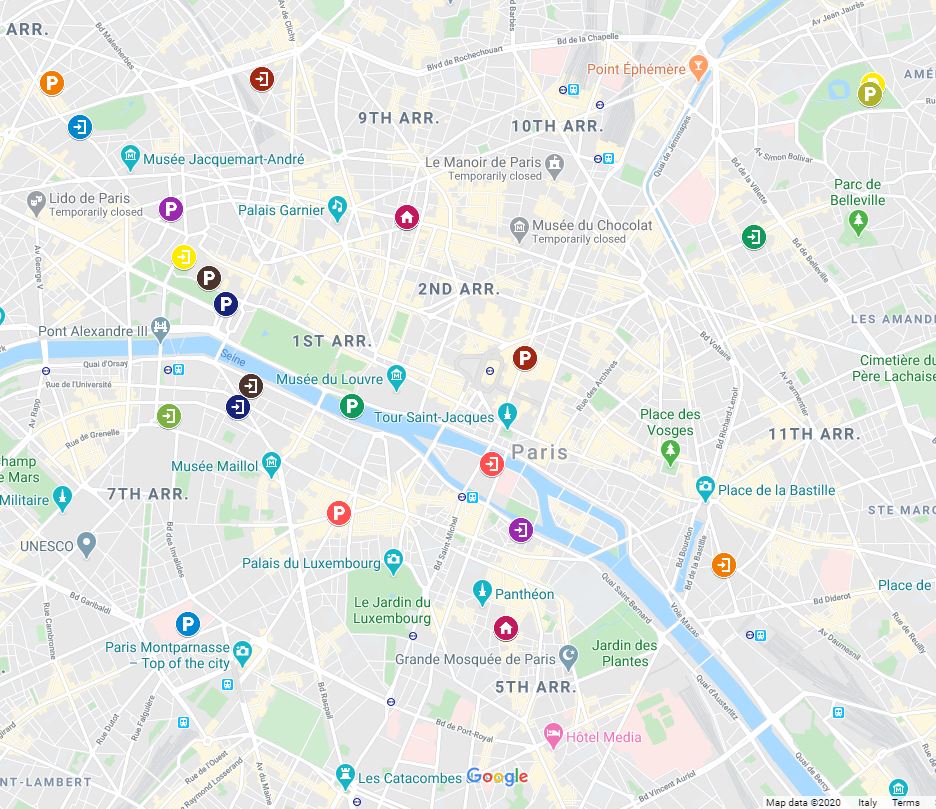}\label{fig:vuota}}
\hspace{2mm}
\subfigure[No collaboration]
{\includegraphics[width=7cm]{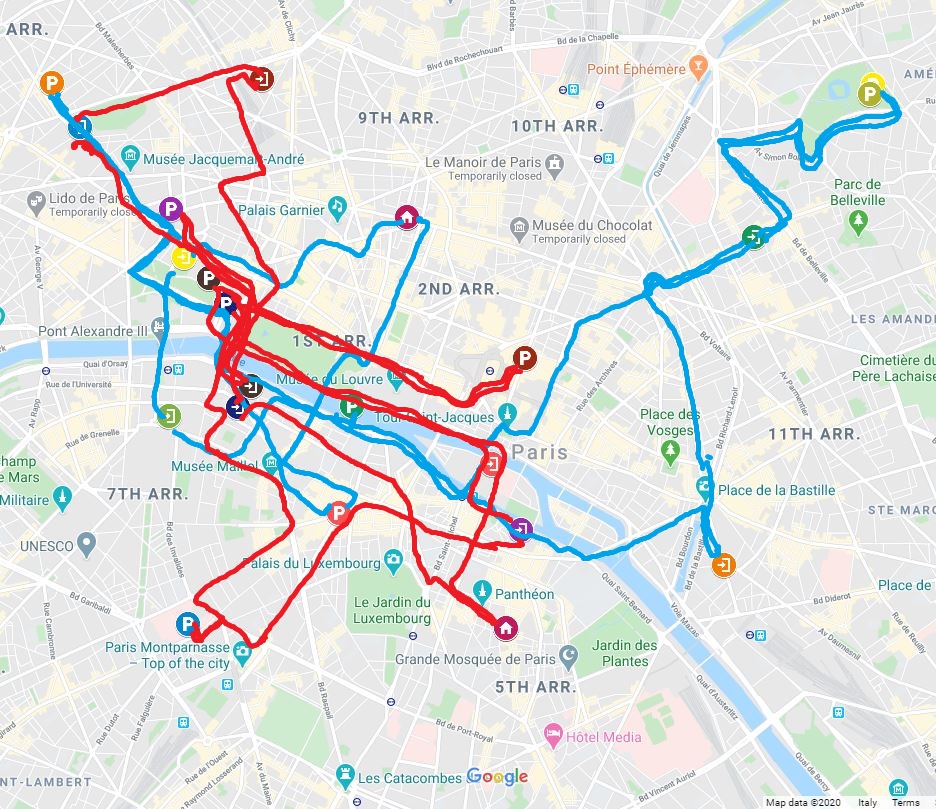}\label{fig:test1}}
\hspace{2mm}
\subfigure[Unconstrained collaboration]
{\includegraphics[width=7cm]{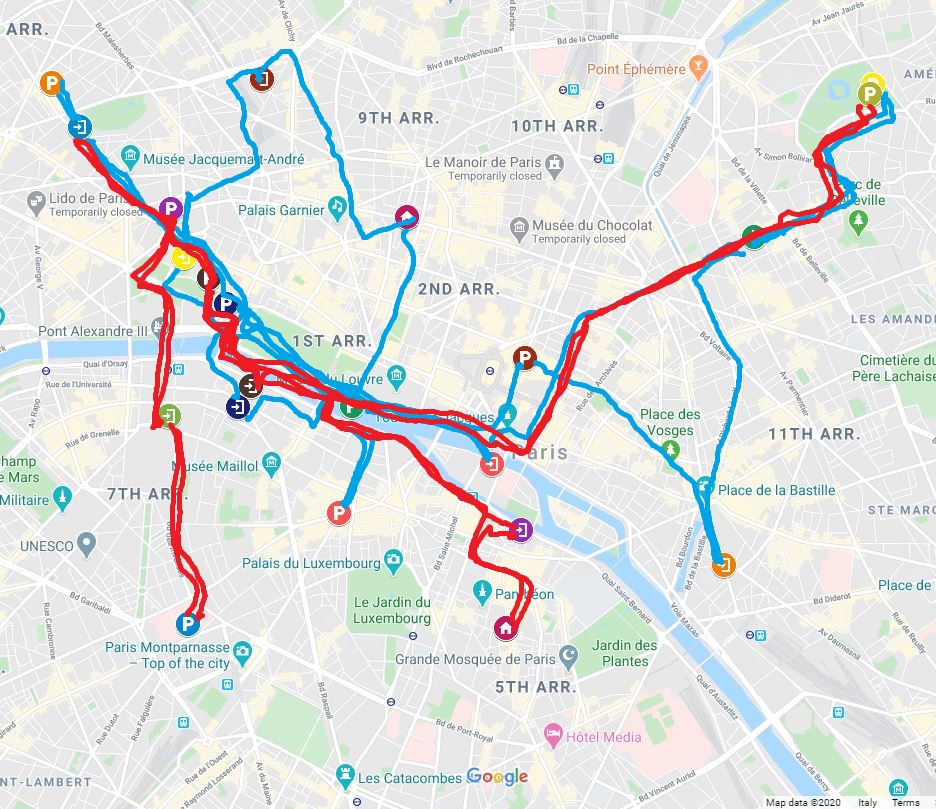}\label{fig:test2}}
\hspace{2mm}
\subfigure[Constrained collaboration]
{\includegraphics[width=7cm]{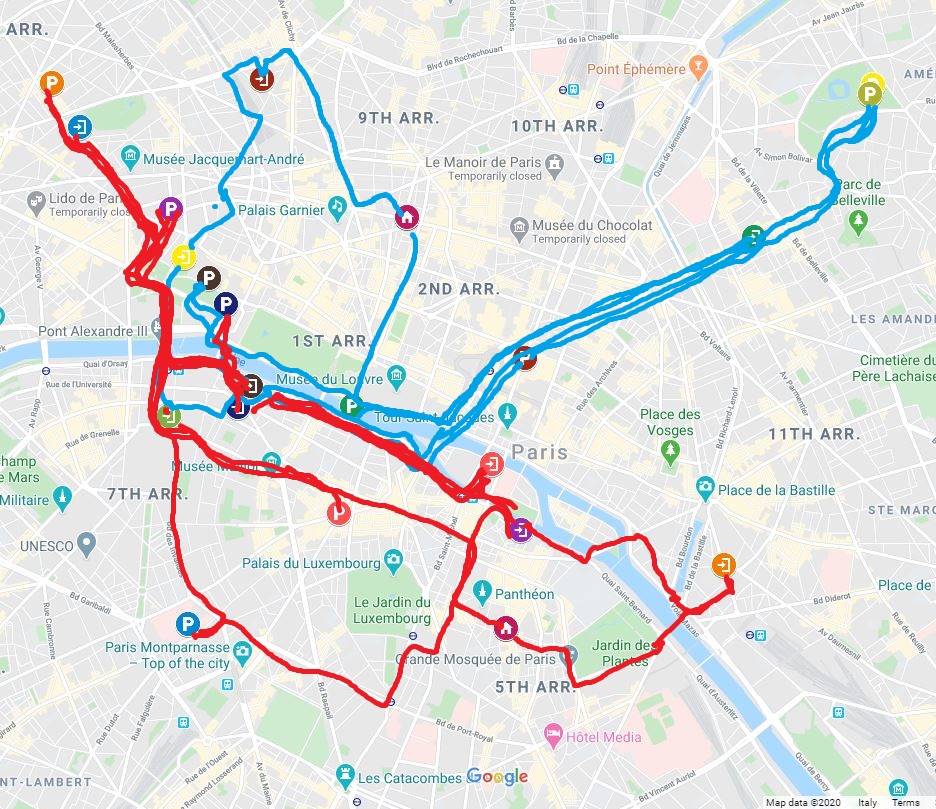}\label{fig:test3}}
\caption{Requests assignment to companies in different operating scenarios 
(pick up \protect\Pickup , delivery points \protect\Delivery, depots \protect\Depot )
}
\end{figure}

In order to investigate the benefits of horizontal collaboration among transportation companies, a map-based example related to the city of Paris is presented with customers and depots located as in Figure \ref{fig:vuota}.

In Figure \ref{fig:test1}, the route plan of two different companies (blue and red) operating a shared taxi service is shown where, as in usual non-collaborative schemes, each company is servicing its own customers with a total routing cost (in minutes) of 98.85.
Conversely, in Figure \ref{fig:test2} the unconstrained route plan, as described in \cite{molenbruch2017benefits} and where each company can serve any other companies' customers, is depicted. In this case the total routing cost drops to 87.5 minutes with a 11.5\% saving on the overall cost. We can observe that exchanging customers may introduce some undesired levels of workload imbalance  in terms of number of served customers of the two companies.
Indeed, two red customers are now serviced by the blue company, and one blue customer is  serviced by the red company, so that, with respect to non-collaborative scenario, the blue company is servicing one customer more and the red company is servicing one customer less.

If we impose that the number of customers served by each company remains the same they had in the non-collaborative setting, the collaboration could remain fruitful for the system and fair to the individual companies as well.
This situation is depicted in Figure \ref{fig:test3} where two red customers are serviced by the blue company, and two blue customers are serviced by the red company with a coalition gain of 10.7\%, just 0.8\% less than in the unconstrained setting.


This example suggests that a balanced approach to collaboration may be more satisfactory for the involved companies and foster the acceptance of the system while, at the same time, provide solutions with costs close to those obtained in an unconstrained setting. The concept shown in the example for the number of customers could  also be shown for the time traveled by the companies to serve other companies' customers.

\section{Models for a fair collaboration}\label{models}

In this section we assume that a coalition has been established among companies that transport one or more customers from origin to destination, and present two models for the optimization of the transportation operations with  constraints aimed at controlling unfairness among  companies.

The first model focuses on balancing, for each company, the time spent in servicing requests acquired from other companies with the time needed to serve requests conceded to other companies.
The second model, analogously to the first, balances the number of exchanged customers.

We first introduce some notation.
Let $M$ be the set of transportation companies in the coalition.
Each company $m \in M$ has the following attributes:
an origin depot $h_m'$ and
a destination depot $h_m''$, which we refer to as $H_m=\left\{h_m',h_m''\right\}$,
a set of vehicles $K_m$,
and a set of customers' requests $C_m$.

We denote
with $H = \bigcup_{m \in M} H_m$ the set of all depots and
with $K = \bigcup_{m \in M} K_m$ the overall coalition fleet,
where each vehicle $k\in K$ has
capacity $Q_k$ and
maximum route duration $T_k$.
Furthermore, let $C = \bigcup_{m \in M} C_m$ the set of all requests in the coalition.
Each request $c \in C$ has the following attributes:
an origin $o_c$ (pick-up point),
a destination $d_c$ (drop-off point),
a demand $p_c$ (number of customers that intend to travel from $o_c$ to $d_c$),
a minimum travel time $t_c$ (direct transfer from $o_c$ to $d_c$),
a service time $s_c'$ and $s_c''$ at origin and destination,
a time window $[e_c',l_c']$ and $[e_c'',l_c'']$ at origin and destination, respectively,
and a maximum time allowed on board  $T_c$ (obviously, $T_c \geq t_c$).
The sets of all origin and destination points are
$O = \left\{o_c\ |\ c \in C\right\}$ and
$D = \left\{d_c\ |\ c \in C\right\}$, respectively. In Table \ref{TBL:PBLNOTATION} we resume the notation used to describe the problem.

\begin{table}[tp]
\resizebox{\textwidth}{!}{
\begin{tabular}{|p{3cm}|p{12cm}|}
\hline
\multicolumn{2}{|l|}{\textbf{Sets}}\\
$M$ & set of companies in the coalition\\
$C_m$ & set of requests of company $m \in M$\\
$C$ & set of all requests, $C = \bigcup_{m \in M} C_m$\\
$H_m=\left\{h'_m,\ h''_m\right\}$ & starting and ending depot, for company $m\in M$\\
$H$ & set of all depots $H = \bigcup_{m \in M} H_m$\\
$K_m$ & set of vehicles of company $m \in M$\\
$K$ & set of all vehicles, $K = \bigcup_{m \in M} K_m$\\
\multicolumn{2}{|l|}{\textbf{Requests attributes}}\\
$o_c$ & origin of request $c \in C$\\
$d_c$ & destination of request $c \in C$\\
$O$ & set of origins $O = \left\{o_c\ |\ c\in C\right\}$\\
$D$ & set of destinations $D = \left\{d_c\ |\ c\in C\right\}$\\
$p_c$ & number of customers in request $c \in C$\\
$t_c$ & minimum travel time (direct transfer) from $o_c$ to $d_c$ of request $c \in C$\\
$s'_c$ & service time at origin $o_c$ for request $c \in C$\\
$s''_c$ & service time at destination $d_c$ for request $c \in C$\\
$T_c$ & maximum time allowed on board  for request $c \in C$\\
$[e'_c,l'_c]$ & service time window at origin $o_c$ for request $c \in C$\\
$[e''_c,l''_c]$ & service time window at destination $d_c$ for request $c \in C$\\
\multicolumn{2}{|l|}{\textbf{Vehicles attributes}}\\
$Q_k$ & capacity of vehicle $k\in K$\\
$T_k$ & maximum route duration for vehicle $k\in K$\\
%
\multicolumn{2}{|l|}{}\\
\hline
\end{tabular}}
\caption{Notation - Part I} \label{TBL:PBLNOTATION}
\end{table}

The problem is modeled on a
complete directed graph $G =(V,A)$ where
$V = H \bigcup O \bigcup D$ and
$|V|=2\cdot|M|+2\cdot|C|$.
Thus, each
node $i \in V$ represents either a depot or an origin or a destination, while
arc $(i,j) \in A$ represents a direct transfer from node $i$ to node $j$ with an associated
traveling time  $t_{ij}$ and cost $c_{ij}$.
In more detail,
each node $i\in V$ representing an origin point for a request $c$, i.e. $i=o_c$ for some $c\in C$, inherits
service time $s_i=s_c'$,
time window $[e_i,l_i]=[e_c',l_c']$, and has an associated
flow quantity $q_i = p_c$;
each node $i\in V$ representing a destination  for a request $c$, i.e. $i=d_c$ for some $c\in C$, inherits
service time $s_i=s_c''$,
time window $[e_i,l_i]=[e_c'',l_c'']$, and has an associated
flow quantity $q_i = -p_c$;
each node $i\in V$ representing a depot for a company $m$, i.e. $i=h_m'$ or $i=h_m''$ for some $m\in M$, has
service time $s_i=0$,
time window $[e_i,l_i]=[0,\infty]$, 
and a flow quantity $q_i = 0$.

The decision variables we need are introduced as follows.
\begin{itemize}
\item
Binary variables $y_c^k$ take value 1 if request $c\in C$ is serviced by vehicle $k\in K$ and 0 otherwise.
\item
Binary variables $x_{ij}^k$ take value 1 if arc $(i,j)\in A$ is traversed by vehicle $k\in K$ and 0 otherwise.
\item
Integer variables $w_{i}^k$ represent the number of onboard customers when vehicle $k\in K$ leaves node $i\in V$. 
\item
Continuous variables $u_{i}^k$ represent the  arrival time of vehicle $k$ at node $i\in V$. 
\item
Continuous variables $r_{c}^k$ represent the on board time of request $c\in C$ (excluding service time) when serviced by vehicle $k\in K$.
\end{itemize}

We first present the part of the formulation which is common to the two optimization models we propose.


\begin{align}
 & \min z= \quad \sum_{i \in V} \sum_{j \in V} \sum_{k \in K} c_{ij}x_{ij}^k \label{MODEL:OF}\\
\nonumber & \text{+ Routing constraints}\\
\nonumber & \text{+ Time windows and maximum time on board  constraints}\\
\nonumber & \text{+ Capacity constraints}\\
&  y_{c}^k \in \{0,1\} & k \in K \quad  c \in C \label{MODEL:20}\\
&  x_{ij}^k \in \{0,1\} &  k \in K \quad i \in V \quad  j \in V \label{MODEL:15}\\
&  w_{i}^k \in \mathbb{N} &  k \in K \quad  i \in V \label{MODEL:16}\\
&  u_{i}^k \ge 0 &  k \in K \quad  i \in V \label{MODEL:18}\\
&  r_{c}^k \ge 0 &  k \in K \quad  c \in C. \label{MODEL:17}
\end{align}

The objective function minimizes the total cost of the transportation operations for the coalition. The detailed explanation of routing, time windows and maximum ride time and capacity constraints follows.

\subparagraph{\textit{Routing constraints}.}

\begin{align}
& \sum_{j \in V} x_{h_m',j}^k = \sum_{i \in V} x_{i,h_m''}^k = 1 &  m \in M \quad  k \in K_m \label{MODEL:3}\\
& \sum_{i \in V} x_{i,l}^k - \sum_{j \in V} x_{l,j}^k = 0 &  l \in V \quad  k \in K \label{MODEL:5}\\
& \sum_{k \in K} y_c^k = 1 &  c \in C\label{MODEL:1}\\
& \sum_{j \in V} x_{o_c,j}^k=y_c^k &  c \in C \quad  k \in K \label{MODEL:13}\\
& \sum_{j \in V} x_{o_c,j}^k - \sum_{i \in V} x_{i,d_c}^k = 0 &  c \in C \quad  k \in K \label{MODEL:4}
\end{align}

Constraints \eqref{MODEL:3} guarantee that each vehicle $k\in K_m$ traverses exactly one arc exiting from the origin depot $h_m'$ and exactly one arc entering the destination depot $h_m''$.
In other words, each vehicle is allowed only one one-way trip between its company's depots. This tour is empty  when variable $x_{h_m',h_m''}^k=1$.
Note that in most cases origin and destination depots are physically the same, though for modelling convenience they are considered as distinct.
%
%
Constraints \eqref{MODEL:5} guarantee flow conservation, that is, if vehicle $k$ enters a node $l$ then it must also leave the same node.
Constraints \eqref{MODEL:1} guarantee that each request $c$ is serviced by exactly one vehicle of the coalition.
Constraints \eqref{MODEL:13} guarantee that vehicle $k$ exits pick-up node $o_c$ if and only if it serves request $c$.
Constraints \eqref{MODEL:4} force vehicle $k$ to visit pick-up node $o_c$ if and only if it also visits drop-off node $d_c$.
Subtours elimination and precedence constraints are implied by the following constraints on traveling times.

\subparagraph{\textit{Time windows and maximum on board time constraints}.}

\begin{align}
& u_j^k \ge (u_i^k + s_i + t_{ij})x_{ij}^k &  i,j \in V \quad  k \in K \label{MODEL:6}\\
& r_c^k = u_{d_c}^k - (u_{o_c}^k + s_{o_c}) &  c \in C \quad  k \in K \label{MODEL:7}\\
& u_{h_m''}^k - u_{h_m'}^k \le T_k &  m \in M \quad  k \in K_m \label{MODEL:8}\\
&  e_{i} \le  u_{i}^k \le  l_{i} &  i \in V \quad  k \in K \label{MODEL:9}\\
&  t_{o_c,d_c} \le r_{c}^k &  c \in C \quad  k \in K \label{MODEL:10bis}\\
&  r_{c}^k \le  T_c &  c \in C \quad  k \in K \label{MODEL:10}
\end{align}

Constraints \eqref{MODEL:6} ensure the consistency of arrival times at nodes visited by vehicle $k$.
Constraints \eqref{MODEL:7} evaluate the on board time for each request $c$ as the difference between the arrival time at drop-off node $d_c$  and the departure time (arrival time plus service time) from pick-up node $o_c$.
Constraints \eqref{MODEL:8} bound the route duration for each vehicle $k$ at a value not greater than its maximum allowed value $T_k$.
Constraints \eqref{MODEL:9} ensure that arrival time at each node $i$ is in the required time window.
Constraints \eqref{MODEL:10bis}  bound the on board time to be at least the direct traveling time $t_{o_c,d_c}$ and forbids the tour to visit the destination before the origin.
Constraints \eqref{MODEL:10} bound the on board time of each request $c$ to be not greater than its maximum allowed value $T_c$.\\

\subparagraph{\textit{Capacity constraints}.}

\begin{align}
& w_j^k \ge (w_i^k + q_j)x_{ij}^k &  i,j \in V \quad  k \in K \label{MODEL:11}\\
&  \max \{0 , q_i \} \le  w_{i}^k \le  \min \{ Q_k,Q_k+q_i\} &  i \in V \quad  k \in K \label{MODEL:12}\\
&  w_{h_m'}^k = 0 &  m \in M \quad  k \in K_m \label{MODEL:12bis}
\end{align}

Constraints \eqref{MODEL:11} ensure the consistency of vehicle load at nodes visited by vehicle $k$.
Constraints \eqref{MODEL:12} ensure that the capacity of each vehicle $k$ is not violated during its trip.
Constraints \eqref{MODEL:12bis} fix the initial load of each vehicle $k$ when leaving the depot.\\

\subparagraph{\textit{Constraints linearization}.}

Constraints (\ref{MODEL:6}) and (\ref{MODEL:11}) are nonlinear. These sets of constraints can be rewritten in linear form in the problem variables using the big-M technique (as suggested in \cite{desrochers1987vehicle}, \cite{desrosiers1995time} and \cite{desrochers1991improvements}) as follows.

\bigskip
\begin{align}
& u_j^k \ge u_i^k + s_i + t_{ij}-U_{ij}^k(1-x_{ij}^k) &  i,j \in V \quad  k \in K \label{MODEL:21}\\
& w_j^k \ge w_i^k + q_i-W_{ij}^k(1-x_{ij}^k) &  i,j \in V \quad  k \in K \label{MODEL:22}
\end{align}
\noindent
where $U_{ij}^k$ and $W_{ij}^k$ parameters have to be chosen according to the following rules: $U_{ij}^k \ge \max \{0,l_i+s_i+t_{ij}-e_j\}$ and $W_{ij}^k \ge \min \{Q_k,Q_k+q_i\}$.
\bigskip





\paragraph{Handling special customers.}
In order to enhance the collaboration between companies it is reasonable to give companies the possibility to a priori serve some of  their own customers or to refuse some of  other companies.
For each request $c \in C \setminus C_m$ refused by company $m$, models can be easily adapted by adding a new constraint
$$\sum_{k \in K_m} y_c^k = 0.$$
On the other hand, if company $m$ wants to serve its request $c \in C_m$, constraint
$$\sum_{k \in K_m} y_c^k = 1$$
can be added.

\subsection{The balanced collaborative DARP models}\label{ES}

In this section we introduce the models aimed at controlling the workload balance in terms of time and/or number of customers.

For the sake of clarity, we say that, from the point of view of company $m\in M$, a request $c\in C$ is \textit{conceded} if $c\in C_m$ and $c$ is served by a different company $m'\neq m$. A request $c\in C$ is  \textit{acquired} if $c$ is served by the company $m$ and $c\in C_{m'}$ for some $m'\neq m$.

\paragraph{The \TCDARP\ model.}

In the \textit{Time balanced Collaborative DARP} model (\TCDARP)  the \textit{time balance} $S_m$ for company $m\in M$ is computed as the difference between the sum of minimum travel times of all acquired requests and the sum of minimum travel times of all conceded requests.
We recall that the minimum travel time $t_c$ of a request $c\in C$  is the travel time of a direct transfer from $o_c$ to $d_c$.

We aim at keeping the absolute value of the time balance $S_m$ below a predefined threshold $\widetilde{S}_m$.
The  choice of the value of $\widetilde{S}_m$ will be
part of the agreement among the companies involved in the coalition.
For instance, it might be set to a fixed percentage $\alpha_T$ of the total minimum travel times of requests owned by company $m$, namely $\widetilde{S}_m=\alpha_T\cdot\sum_{c \in C_m} t_c $.

 We also consider a term $S_m'$ which represents an offset for the time balance. It may be seen as a past credit, if positive (debt, if negative), that companies have agreed to consider in the balance. We say that the time balance is \textit{without memory} when $S_m'$ is set to 0, \textit{with memory} otherwise.

\medskip
The \TCDARP\ model is, thus, defined as follows.\\

\noindent\textit{\TCDARP\ model}
\begin{align}
\nonumber & \Theta^*_{\TCDARP} = \min \quad \sum_{i \in V} \sum_{j \in V} \sum_{k \in K} c_{ij}x_{ij}^k\\
\nonumber & \text{+ Routing constraints}\\
\nonumber & \text{+ Time windows and maximum ride time constraints}\\
\nonumber & \text{+ Capacity constraints}\\
\nonumber & \text{+ Decision variables domain constraints \eqref{MODEL:20}-\eqref{MODEL:17}}\\
& S_m = \sum_{k \in K_m} \sum_{c \in C\setminus C_m} t_cy_c^k-\sum_{k \in K\setminus K_m} \sum_{c \in C_m} t_cy_c^k &  m \in M\label{MODEL:33bis}\\
& -\widetilde{S}_m \le S_m + S_m' \le \widetilde{S}_m &  m \in M \label{MODEL:33}
\end{align}

Constraints \eqref{MODEL:33bis} compute the time balance $S_m$ for each company $m\in M$.

Constraints \eqref{MODEL:33} bound the absolute value of the time balance (including the possible offset $S_m'$) of each company $m$ to be lower than a fixed threshold $\widetilde{S}_m$ which may differ for each company.

\paragraph{The \CCDARP\ model.}

An alternative way to balance the activities of the companies involved in the coalition is to balance the requests served. The \textit{Customer  balanced Collaborative DARP} model (\CCDARP)  aims at bounding the balance of customers served by each company.
The \textit{customer balance} $U_m$ is computed as the difference between the total number of acquired requests and that of conceded requests. Its absolute value is constrained to be not greater than a predefined threshold $\widetilde{U}_m$.

Analogously to the \TCDARP\ model,
the value of $\widetilde{U}_m$ might be set as a fixed percentage $\alpha_C$ of the total number of customers owned by company $m$, namely $\widetilde{U}_m=\alpha_C\cdot\sum_{c \in C_m} p_cy_c^k$;
accordingly, we also introduce an offset  $U_m'$ for the customer balance $U_m$ and say that the model is \textit{without memory} when $\bar{U_m}$ is set to 0 and \textit{with memory} otherwise.

\medskip
The model can be, thus, formulated replacing constraints (\ref{MODEL:33bis},\ref{MODEL:33}) with constraints related to the customers.\\

\noindent\textit{\CCDARP\ model}

\begin{align}
\nonumber & \Theta^*_{\CCDARP} = \min \quad \sum_{i \in V} \sum_{j \in V} \sum_{k \in K} c_{ij}x_{ij}^k\\
\nonumber & \text{+ Routing constraints}\\
\nonumber & \text{+ Time windows and maximum ride time constraints}\\
\nonumber & \text{+ Capacity constraints}\\
\nonumber & \text{+ Decision variables domain constraints \eqref{MODEL:20}-\eqref{MODEL:17}}\\
& U_m =\sum_{k \in K_m} \sum_{c \in C\setminus C_m} p_cy_c^k-\sum_{k \in K\setminus K_m} \sum_{c \in C_m} p_cy_c^k &  m \in M \label{MODEL:35bis}\\
& -\widetilde{U}_m \le  U_m + U_m' \le \widetilde{U}_m \hspace{1.7cm} & m \in M \label{MODEL:35}
\end{align}

Constraints (\ref{MODEL:35bis})  compute the customer balance while constraints (\ref{MODEL:35}) bound its absolute value.

\medskip
It may be also of interest to add constraints (\ref{MODEL:35bis},\ref{MODEL:35}) to the \TCDARP\ model instead of substituting (\ref{MODEL:33bis},\ref{MODEL:33}). In this case the model will be denoted as \TCCDARP.
In Table \ref{csoh:notation} we summarize the notation used to formulate the models.

\begin{table}[tp]
\resizebox{\textwidth}{!}{
\begin{tabular}{|p{3cm}|p{12cm}|}
\hline
\multicolumn{2}{|l|}{\textbf{Model parameters}}\\
$V = H\bigcup O \bigcup D$ & set of vertices\\
$A$ & set of arcs\\
$t_{ij}$ & traveling time on arc $\left(i,j\right) \in A$\\
$c_{ij}$ & cost of arc $\left(i,j\right) \in A$\\
$s_i$ & service time in node $i \in V$\\
$[e_i,l_i]$ & time windows at node $i \in V$\\
$q_i$ & flow of customers at node $i \in V$\\
$S_m'$ & time balance offset for company $m\in M$\\
$U_m'$ & customers balance offset for company $m\in M$\\
\multicolumn{2}{|l|}{\textbf{Companies attributes}}\\
$\widetilde{S}_m$ & maximum time balance allowed for company $m\in M$\\
$\widetilde{U}_m$ & maximum customers balance allowed for company $m\in M$\\
$\alpha_T$ & parameter suggested to define $\widetilde{S}_m$ for company $m\in M$\\
$\alpha_C$ & parameter suggested to define $\widetilde{U}_m$ for company $m\in M$\\
%
\multicolumn{2}{|l|}{\textbf{Decision variables}}\\
$y_c^k$ & 1 if request $c \in C$ is performed by vehicle $k\in K$, 0 otherwise\\
$x_{ij}^k$ & 1 if arc $(i,j)\in A$ is used by vehicle $k\in K$, 0 otherwise\\
$w_{i}^k$ & customers on vehicle $k\in K$ after visiting node $i\in V$\\
$u_{i}^k$ & arrival time of vehicle $k\in K$ at node $i\in V$\\
$r_{c}^k$ & on board time of request $c \in C$ serviced by vehicle $k\in K$\\
$S_m$ & time balance for company $m\in M$\\
$U_m$ & customers balance for company $m\in M$\\
\multicolumn{2}{|l|}{}\\
\hline
\end{tabular}}
\caption{Notations - Part II} \label{csoh:notation}
\end{table}

\section{An Adaptive Large Neighborhood Search algorithm}\label{heur}

The Adaptive Large Neighborhood Search (ALNS) algorithmic approach has been widely used in the literature, especially for routing and scheduling problems (see \cite{muller2009adaptive} and \cite{pisinger2010large} and references therein). The ALNS scheme, first proposed in \cite{ropke2006adaptive}, is an extension of the Variable Neighbourhood Search (VNS, see \cite{vns} for details) in which a number of destroy and repair operators are iteratively and pseudo-randomly applied in order to improve the current solution.
The ALNS keeps track of the frequency of success of each destroy and repair operator and updates a parameter which is used as probability to choose the operator in the next  iteration.
In order to implement an ALNS heuristic, one out of many metaheuristic frameworks can be chosen (simulated annealing, tabu search, guided local search, etc.).
We chose the simulated annealing since it is known to be successful in routing problems, as assessed by \cite{ropke2006adaptive}.

The ALNS heuristic is outlined in Algorithm \ref{alsn}.
The algorithm takes in input a feasible solution used to initialize both the current solution $x$ and the incumbent solution $x^*$ (best current solution).
At each iteration, a destroy operator and a repair operator are drawn from their respective pools and applied in sequence to partially destroy solution $x$ and possibly construct a new feasible solution $x'$. If the construction phase fails, the original solution $x$ is returned.
Then, the current solution $x$ may be replaced by $x'$ according to a probability that depends on how far the cost of the new solution $x'$ is  from the cost of the incumbent $x^*$ and a temperature parameter $T_{max}$ which is lowered at each iteration according to a multiplier $\gamma$.
The heuristic continues to iterate until $T_{max}$ becomes less than or equal to 1, then the ALNS stops.
According to the used rule, if solution $x'$ improves upon the incumbent, then the probability is automatically set to 1, while if the solution $x'$ is very bad, its cost is very large and the acceptance probability is very small.

Destroy and repair operators are chosen independently with a roulette-wheel selection mechanism, as explained in \cite{ropke2006adaptive}.
Namely, destroy and repair operators are drawn from their respective pools according to a score that is initialized to be equal for each of them, and later updated and increased each time an operator successfully contributes to find an improving solution.
%
%
The scores need to be reset during the process to escape from situations in which a destroy/repair operator, that was successful in the first steps, gains a very high score and prevents other operators from being selected even though the operator itself may stop being effective and should leave more opportunities to others.
Every time an improving solution is found, a counter $r$ is increased and, when the counter $r$ reaches a refresh threshold $R$, all scores are reset to the same value and destroy/repair operators become again equally probable.

According to \cite{pisinger2010large}, 
the degree of destruction of destroy operators is quite crucial since a too small degree may result in restricting the exploration of the search space,
while a too large degree easily degrades into repeated re-optimizations.
In \cite{shaw1998using}, it is suggested to gradually increase the degree of destruction while, in \cite{ropke2006adaptive}, it is suggested to choose, at each iteration, the degree of destruction randomly in a given interval.
The latter mechanism is also used in \cite{molenbruch2017benefits}.
In Algorithm \ref{alsn} we implement an adaptive mechanism to  choose the degree of destruction $q$.
As suggested by \cite{ropke2006adaptive}, a minimum and a maximum number of requests, namely $q_{min}$ and $q_{max}$, that can be removed at each iteration are provided as input parameters to the algorithm.
The heuristic starts by setting $q$ at the minimum value $q_{min}$.
Then, for each iteration in which the heuristic is not able to find an improving solution, a counter $w$ is incremented until an enlarging threshold $E$ is reached.
When the enlarging threshold is passed, the degree of destruction $q$ is incremented unless it is already at its maximum value.
In order to avoid  increasing $q$ too much, at each iteration the value might be reduced according to a small fixed probability $p$.
Counter $w$ is reset whenever $q$ is modified.
The procedure for resizing the parameter $q$ is outlined in Algorithm \ref{q}.

\normalem
\begin{algorithm}
	\caption{ALNS scheme}
	\label{alsn} {
		\DontPrintSemicolon
		\SetKwInOut{Input}{input}\SetKwInOut{Output}{output}\SetKwInOut{Global}{global}
		\Input{$x_0$,$T_{max}$,$\gamma$,$R$,$q_{min}$,$q_{max}$,$p$,$E$}
		\Output{$x^*$ heuristic solution}

        \medskip

		-- $x=x_0$ as current solution and $x^*=x_0$ incumbent solution;\\
		-- $q:=q_{min}$;\\
		-- Initialize $r$ and $w$ to zero;\\
		-- Set destroy and repair operators' scores to initial values;\\
        \medskip
		\While{$T_{max}>1$}{
					-- $q:=resizeNeighborhood(E,w,q,q_{min},q_{max},p)$;\\
					-- Draw a destroy and a repair operator;\\
					-- Destroy the current solution $x$; \\
					-- Repair the destroyed solution and obtain $x'$;\\
					-- \If{ $cost(x')\ge cost(x^*)$}{--$w=w+1$;}
					-- $u$ randomly drawn in $U(0,1)$;\\
					\If{$u<e^{\frac{cost(x^*)-cost(x')}{T_{max}}}$}{
					-- $x=x'$;\\
						\If{$cost(x')<cost(x^*)$}{
						-- $x^*=x'$;
						
						-- \If{$r>R$}{
						-- $r=0$;\\
						-- Set destroy and repair operators' scores to initial values;\\
					}
					\Else{
						-- Update destroy and repair operators' scores;\\}
						-- $r=r+1$;
						}
						}
					
					-- $T_{max}=T_{max}*\gamma$;\\

					}
    	-- \textbf{return} $x^*$;
    }
\end{algorithm}

\normalem
\begin{algorithm}
	\caption{$resizeNeighborhood$ function}
	\label{q} {
		\DontPrintSemicolon
		\SetKwInOut{Input}{input}\SetKwInOut{Output}{output}\SetKwInOut{Global}{global}
		\Input{$E,w,q,q_{min},q_{max},p$}
		\Output{($q,w$) }

        \medskip

    	-- \If{$w>E$ and $q<q_{max}$}{-- $q=q+1$;\\--$w=0$;}
					-- $reduc$ randomly drawn in $U(0,1)$;\\
					-- \If{$reduc<p$ and $q>q_{min}$}{-- $q=q-1$;\\--$w=0$;\\}
					-- \textbf{return} $q$;
    }
\end{algorithm}

\subsection{Destroy and repair operators}

In the proposed ALNS heuristic six destroy and three repair operators have been identified. One repair operator is parametric, and we used three different parameter values to obtain a total of five repair operators.

\subsubsection{Destroy operators}

The following destroy operators have been used: the random removal, the worst removal, the related removal (also known as Shaw removal), the proximity removal, the closeness removal, and the interchangeability removal.
The former three destroy operators are well known in the literature and come from the seminal work of \cite{ropke2006adaptive}.
The fourth operator was introduced by \cite{molenbruch2017benefits} whereas the fifth and sixth destroy operators are a new contribution of this paper.
Every operator takes as input the current solution and parameter $q$ that defines the degree of destruction.

\begin{itemize}
\item \textbf{\textit{Random removal.}}
    All requests are assigned the same probability $p=1/|C|$ and $q$ distinct requests are randomly selected according to this distribution. Selected requests are removed from the solution.

\item \textbf{\textit{Worst removal.}}
    For each request $a \in C$, a marginal cost $mc_a$ is computed as the difference between the cost of the solution servicing all requests and the cost of the solution obtained by removing request $a$.
    Then, a probability $p_a = mc_a^2/(\sum_{c \in C}mc_c^2)$ is assigned to each request $a \in C$, and $q$ distinct requests are randomly selected and removed from the solution. Requests with higher marginal cost have more probability to be selected and removed.

\item \textbf{\textit{Related removal.}}
    A request $r\in C$ is first selected with uniform probability $p=1/|C|$.
    Then, a relatedness measure (with respect to $r$)
    $$rel_{ar}=\left(\frac{t_{o_a,o_r}+t_{d_a,d_r}}{\max_{(ij) \in A} t_{ij}}+\frac{|\frac{l_{o_a}-e_{o_a}}{2}-\frac{l_{o_r}-e_{o_r}}{2}|+|\frac{l_{d_a}-e_{d_a}}{2}-\frac{l_{d_r}-e_{d_r}}{2}|}{\max_{(ij) \in A} l_i-\min_{(ij) \in A}e_i}\right)^{-1},$$%
    is assigned to each remaining request $a\in C\setminus \{r\}$.
    Finally, a probability $p_a = rel_{ar}^2/(\sum_{c \in C\setminus \{r\}}rel_{cr}^2)$
    is assigned to each request $a\in C\setminus \{r\}$, and $q-1$ more distinct requests are randomly selected. The $q$ selected requests are then removed from the solution.

    The relatedness measure considers both spatial and temporal similarity between two requests. According to \cite{ropke2006adaptive}, extracting two related requests from a feasible solution allows more chance to rebuild a new feasible solution during the reconstruction phase.
    Note that the relatedness measure between any pair of requests does not depend on the current solution and, thus, can be computed a priori before starting the iterative part of the ALNS heuristic.

\item \textbf{\textit{Proximity removal.}}
    For each request $a \in C$, a proximity measure $prox_a = \min_{c \in D}\left(rel_{a,c}\right)$ is computed, where $D$
    is the set of requests on a different route with respect to $a$ and $rel_{a,c}$ is the relatedness measure explained above.
    Then, a weight $p_a = prox_a^2/(\sum_{c \in C\setminus \{a\}}prox_c^2)$ is assigned to each request $a\in C$ and $q$ distinct requests are randomly selected according to these weights and removed from the solution.

    Proximity of a request $a$ in a solution is built upon relatedness between requests and considers the minimum relatedness with requests in different routes. This should enhance horizontal collaboration by extracting requests with high opportunity to be inserted at low cost in other routes during the reconstruction phase.

\item \textbf{\textit{Closeness removal.}}
   This is a new operator we introduce as a variant of related removal where the relatedness measure between pairs of requests is replaced by a closeness measure.
    The closeness measure $close_{ar}$ of request $a \in C$ to request $r \in C$ aims at indicating how much request $a$ is compatible with request $r$.
    To compute $close_{ar}$ we first consider the  six possible vertex sequences $S_{ar}$ to serve both and only requests $a$ and $r$.
    Figure \ref{sequences} depicts the six possible visiting sequences in $S_{ar}$.

\vspace{12pt} \noindent
\begin{figure}[htbp!]
\begin{center}
\begin{tikzpicture}[->,>=stealth',shorten >=1pt,auto,node distance=1.5 cm,
thick,main node/.style={circle,fill=blue!20,draw, minimum width=6pt,inner sep=3pt},main2 node/.style={circle,fill=green!20,draw, minimum width=6pt,inner sep=3pt},main3 node/.style={circle,fill=red!20,draw, minimum width=6pt,inner sep=3pt},main4 node/.style={circle,fill=orange!20,draw, minimum width=6pt,inner sep=3pt},main5 node/.style={circle,fill=gray!20,draw, minimum width=6pt,inner sep=3pt},main6 node/.style={circle,fill=yellow!20,draw, minimum width=6pt,inner sep=3pt}]

\node[main node] (0) [] {$O_r$};
\node[main node] (1) [right of=0] {$D_r$};
\node[main node] (2) [right of=1] {$O_a$};
\node[main node] (3) [right of=2] {$D_a$};

\node[main2 node] (4) [right of=3] {$O_a$};
\node[main2 node] (5) [right of=4] {$D_a$};
\node[main2 node] (6) [right of=5] {$O_r$};
\node[main2 node] (7) [right of=6] {$D_r$};

\node[main3 node] (8) [below of=0] {$O_r$};
\node[main3 node] (9) [right of=8] {$O_a$};
\node[main3 node] (10) [right of=9] {$D_r$};
\node[main3 node] (11) [right of=10] {$D_a$};

\node[main4 node] (12) [right of=11] {$O_a$};
\node[main4 node] (13) [right of=12] {$O_r$};
\node[main4 node] (14) [right of=13] {$D_a$};
\node[main4 node] (15) [right of=14] {$D_r$};

\node[main5 node] (16) [below of=8] {$O_r$};
\node[main5 node] (17) [right of=16] {$O_a$};
\node[main5 node] (18) [right of=17] {$D_a$};
\node[main5 node] (19) [right of=18] {$D_r$};

\node[main6 node] (20) [right of=19] {$O_a$};
\node[main6 node] (21) [right of=20] {$O_r$};
\node[main6 node] (22) [right of=21] {$D_r$};
\node[main6 node] (23) [right of=22] {$D_a$};

\path[edge_style/.style={draw=black, ultra thick},every node/.style={font=\sffamily\Large}]

(0) edge node [above,sloped]{} (1)
(1) edge node [above,sloped]{} (2)
(2) edge node [above,sloped]{} (3)
(4) edge node [above,sloped]{} (5)
(5) edge node [above,sloped]{} (6)
(6) edge node [above,sloped]{} (7)
(8) edge node [above,sloped]{} (9)
(9) edge node [above,sloped]{} (10)
(10) edge node [above,sloped]{} (11)
(12) edge node [above,sloped]{} (13)
(13) edge node [above,sloped]{} (14)
(14) edge node [above,sloped]{} (15)
(16) edge node [above,sloped]{} (17)
(17) edge node [above,sloped]{} (18)
(18) edge node [above,sloped]{} (19)
(20) edge node [above,sloped]{} (21)
(21) edge node [above,sloped]{} (22)
(22) edge node [above,sloped]{} (23)
;

\end{tikzpicture}
\end{center}
\caption{The six possible sequences to serve requests $a$ and $r$}
\label{sequences}
\end{figure}
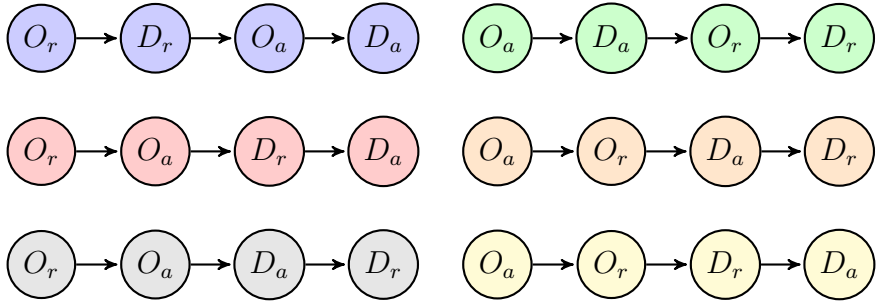

For each sequence $s\in S_{ar}$ which is feasible with respect to time constraints (time windows and travels duration) we compute $close_s=time_s - time_r$, where $time_s$ is the minimum time required to perform all operations in the given order, and $time_r$ is the minimum time required to serve request $r$ only ($t_r$ according to the notation adopted in this paper).
If a sequence $s$ is not feasible we set $close_s=\infty$.
Then, the closeness measure $close_{ar}$ is evaluated as the minimum value among all $close_s$, i.e. $close_{ar}=\min\limits_{s \in S_{ar}} \left(close_s\right) = \min\limits_{s \in S_{ar}} \left(time_s - time_r\right)$.
In other words, the closeness measure $close_{ar}$ provides a lower bound on the time cost of adding request $a$ to a tour  $r$ is already inserted into.
The operator aims at generalizing the relatedness operator by considering a request $a$ close to a request $r$ if it has a low impact on a tour already including $r$.
The Closeness removal operator works exactly like the Related removal, but $close_{ar}$ is used in place of $rel_{ar}$ after the selection of the first request $r$.


Analogously to relatedness, closeness can be computed a priori as it does not depend on the current solution. However, while the relatedness measure can be asymmetric because of a possible asymmetry in travel times, the closeness measure is definitely not symmetric ($close_{ar} \neq close_{ra}$).
%


\item \textbf{\textit{Interchangeability removal.}}
   For each request $a \in C$, an interchangeability measure $int_a = \min_{c \in D}\left(close_{a,c}\right)$ is computed, where $D$
    is the set of requests on a different route with respect to $a$ and $close_{a,c}$ is the closeness measure explained above.
    Then, a weight $p_a = int_a^2/(\sum_{c \in C\setminus \{a\}}int_c^2)$ is assigned to each request $a\in C$ and $q$ distinct requests are randomly selected according to these weights and removed from the solution.

    The interchangeability of a request $a$ in a solution is built upon the closeness between requests and considers the minimum closeness with requests in different routes. As for the proximity removal, this should enhance horizontal collaboration by extracting requests which are more likely to be inserted at low cost in other routes during the reconstruction phase.
\end{itemize}

\subsection{Repair operators}
Given a destroyed solution, where $|C|-q$ requests are feasibly serviced and $q$ requests are excluded, repair operators add the excluded requests to routes in an attempt to create a new solution in which all requests are feasibly serviced.
Adopted operators are the classical operators proposed in \cite{ropke2006adaptive} and \cite{molenbruch2017benefits} along with a new operator based on the closeness measure.
All operators work by assigning a selection probability to requests yet to be inserted in the partial solution. After random selection, the chosen request is inserted at minimum cost: best feasible insertion with respect to solution cost increment.
The process is iterated until all requests are inserted and a new feasible solution is obtained. If an iteration fails to feasibly insert the selected request, the whole process fails and the original solution is returned. In the following we report the peculiarities of each repair operator.

\begin{itemize}
\item \textbf{\textit{Best insertion.}}
    All requests to be inserted in the solution are given the same probability.
\item \textbf{\textit{$k^{th}$-regret insertion.}}
    Each request $a$ to be inserted in the solution is assigned a probability proportional to its $k^{th}$-regret value computed as $r_a^k= \sum_{h=2}^k \left(c_a^h-c_a^1\right)$, where $c_a^1$ is the minimum insertion cost of request $a$ and $c_a^h$ is its $h^{th}$ insertion cost, for $h=2,\ldots,k$.
    %
    %
    %
    %
    %
    %
    %
    We used $k=2,3$ and $4$.

\item \textbf{\textit{Closeness insertion.}}
Each request to be inserted in the solution is assigned a probability inversely proportional to its minimum closeness measure with respect to all requests already in the solution.
    The idea is to first insert  the requests that are more likely to be accommodated in any destroyed route.
\end{itemize}

Regarding the feasibility check of the repaired solutions obtained with repair operators, the one proposed in \cite{molenbruch2017benefits} has been used to assess the feasibility of the trip time for each request and of the vehicles. Then, an ad-hoc feasibility check has been developed to check if time balance and/or customer balance constraints are violated.

\section{Computational results}\label{result}


In this section we present a computational study to assess the effectiveness of the proposed models and the performance of the ALNS heuristic.

Effectiveness of the proposed models is evaluated by comparison with a non-collaborative setting, where companies do not share their customers, and an unconstrained setting, where companies agree to act like a single company minimizing the total cost without any fairness constraint.
The unconstrained setting, which we call \CDARP, is obtained by the linear programming model
\eqref{MODEL:OF}-\eqref{MODEL:17}, \eqref{MODEL:3}-\eqref{MODEL:4}, \eqref{MODEL:7}-\eqref{MODEL:10}, \eqref{MODEL:12}-\eqref{MODEL:12bis}, \eqref{MODEL:21}-\eqref{MODEL:22}. In other words, we remove constraints \eqref{MODEL:33bis} and \eqref{MODEL:33} from \TCDARP\ model.
This MILP provides a lower bound for the collaborative models.
The non-collaborative setting, which we call \NCDARP, is obtained from \CDARP\ by considering one company at the time.
By summing the non-collaborative costs of all companies we obtain an upper bound for the collaborative models.

Collected statistics in the computational experiments are defined as follows.
\begin{itemize}
    \item \SAV: objective function saving w.r.t. the non-collaborative setting.
    \item \GAP: heuristic optimality gap.
    \item $\bar{S}=\frac{1}{|M|}\sum\limits_{m \in M} |S_m|$: average absolute value of the time balance.
    \item $\bar{U}=\frac{1}{|M|}\sum\limits_{m \in M} |U_m|$: average absolute value of the customers balance.
    \item $\widehat{S}=\max\limits_{m \in M} |S_m|$: maximum absolute value of the time balance.
    \item $\widehat{U}=\max\limits_{m \in M} |U_m|$: maximum absolute value of the customers balance.
    \item For each destroy and repair operator, we account for the number of times it has been successful in finding an improving solution.
\end{itemize}

We generated 112 instances grouped in 4 sets of 28 map-based instances each. Instances in the same set have the same number of companies and of requests.
All instances are available at \url{http://or-brescia.unibs.it/instances} and details on their features
are provided in Section \ref{map}.

Section \ref{cvsnc} is dedicated to compute the savings obtained by the collaborative models with respect to the non-collaborative and the unconstrained settings.
%
We also analyze the impact of a collaborative setting over a sequence of consecutive days.

The quality of the solutions obtained by the ALNS heuristic against the optimal solutions is assessed in Section \ref{oneinstance} on the instances of groups A and B.
Section \ref{allinstances} is devoted to summarize the results on the performance of the ALNS heuristic on the larger instances of groups C and D. Results presentation relies on 
graphics.

In Section \ref{operatoreffectiveness} we show the effectiveness of the destroy/reapir operators used in the ALNS heuristic.

Exact solutions of the MILP models were obtained by using CPLEX 12.6.0 on a Windows 64-bit computer with Intel Xeon processor E5-1650, 3.50 GHz, and 16 GB RAM.
%
The same machine was used to run the ALNS heuristic.

\subsection{Map-based random instances}\label{map}

As the problem we study in this paper is new to the literature, no benchmark instances are available for the computational study.

A first set of 48 instances, with  16 to 96 requests, was proposed in \cite{cordeau2006branch} and in \cite{ropke2007models}  for the  single-depot single-company DARP where each request involves only one customer.
These instances were extended in \cite{braekers2014exact} for the multi-depot variant in order to mimic the presence of multiple companies.
However, requests were not associated with a specific company as only the unconstrained setting was considered.
This looks unrealistic in a competitive setting where companies accept to cooperate in view of their own performance improvement and not just for the sake of the system itself.
Thus, we decided to generate new benchmark instances and, in order to have a more realistic setting, we created map-based random instances by extracting real travel times using the framework Graphhopper (repository in \cite{GH}).

The instance generator allows us to generate different instances in terms of size and locations.
The size is measured in number of companies (with one vehicle each) and requests involved. Four different sizes have been considered:\\
-- Group A: 2 companies with 4 requests each;\\
-- Group B: 2 companies with 5 requests each;\\
-- Group C: 4 companies with 12 requests each;\\
-- Group D: 10 companies with 10 requests each.


 For each group, 4 subgroups of 7 instances each have been generated. The 4 subgroups correspond to big European metropolitan areas (Paris, Berlin, London and Rome) where we randomly chose depot addresses of the companies. Note that the subset of instances of Group B related to the city of Rome are also denoted as Group B-R7 as they are further used to mimic and test the models' behaviour on a 7-days scenario.
For each pair of size and city we randomly generated 7 demand scenarios by randomly choosing pick-up and drop-off addresses from the same road network.
%
As requests time windows are concerned, we randomly selected either origin or destination and randomly assigned a time window with a fixed width of about 33 minutes (2000 seconds), while the other time window was left at $\left[0,+\infty\right)$.
Traveling times were taken using the Graphhopper API on the selected network, and for sake of simplicity we assumed the cost $c_{ij}$ to be equal to the travel time $t_{ij}$ for each arc $(i,j)\in A$.

Other parameters were fixed as follows for all instances. Vehicle capacity $Q_k$ and maximum route duration $T_k$ were fixed at 3 and about 5 hours and a half (20000 seconds), respectively, for all $k\in K$.
The number of customers, the maximum time allowed on board, service time at pick-up and drop-off points were fixed to 1, 50 minutes, and 2 minutes, respectively, for all requests $c\in C$.

Thus, we  obtained 4 groups (A,B,C,D) of instances homogeneous in size, with 28 instances each (4 cities, 7 demand scenarios) for a total of 112 instances.

\subsection{Collaborative  vs. non-collaborative settings}\label{cvsnc}
In this section we compute the savings obtained by the collaborative models with respect to the non-collaborative setting and the upper bound on savings given by the unconstrained setting.
%
To this aim an exact solver was used to solve models \NCDARP, \CDARP, \TCDARP, \CCDARP\ and \TCCDARP\ on instances of groups A e B.
%

\begin{figure}[htbp!]
\centering
\begin{tikzpicture}
\begin{axis} [width=8cm,
              symbolic x coords={10\%,20\%,30\%},xtick=data,
legend style={at={(0.5,1.15)},
	anchor=north,legend columns=-1},ymin=0,
    xlabel={$\alpha$},ylabel={\SAV\ (\%)}]
\addplot  plot coordinates{

             (10\%,14.45)
             (20\%,14.45)
             (30\%,14.45)
          };
           \addplot  plot coordinates{

             (10\%,11.20)
             (20\%,12.85)
             (30\%,13.28)
          };
          \addplot  plot coordinates{

             (10\%,11.89)
             (20\%,11.89)
             (30\%,13.96)
          };

          \addplot  plot coordinates{

             (10\%,11.20)
             (20\%,11.89)
             (30\%,13.28)
          };
          \legend{\CDARP,\TCDARP,\CCDARP,\TCCDARP}
\end{axis}
\end{tikzpicture}
\caption{Savings from collaboration  on group A instances}
\label{fig1}
\end{figure}
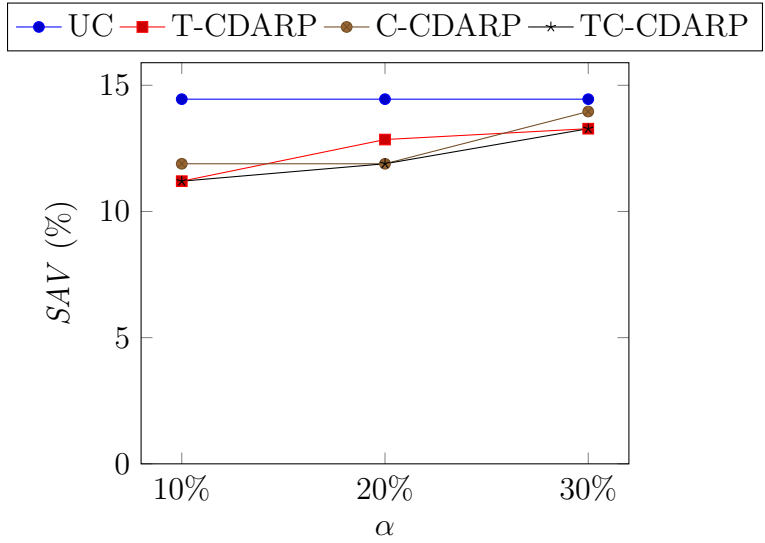

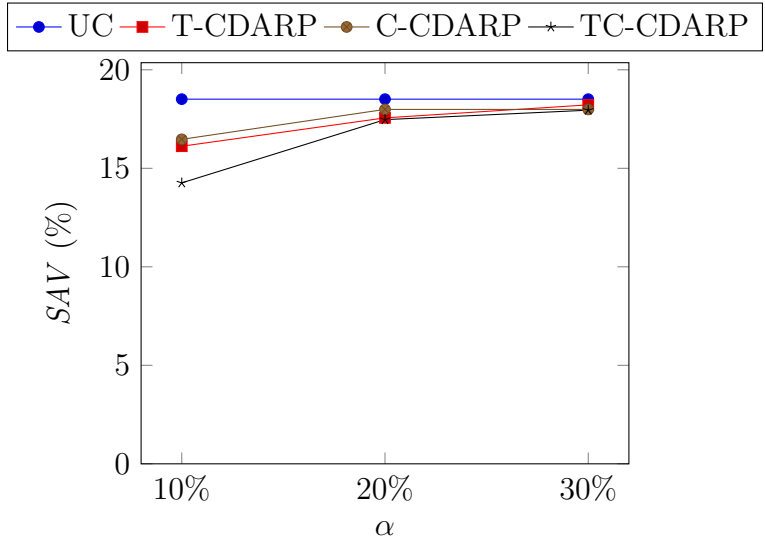
\begin{figure}[htbp!]
\centering
\begin{tikzpicture}
\begin{axis} [width=8cm,
              symbolic x coords={10\%,20\%,30\%},xtick=data,
legend style={at={(0.5,1.15)},ymin=0,
	anchor=north,legend columns=-1},
    xlabel={$\alpha$},ylabel={\SAV\ (\%)}]
\addplot  plot coordinates{

             (10\%,18.51)
             (20\%,18.51)
             (30\%,18.51)
          };
          \addplot  plot coordinates{

             (10\%,16.12)
             (20\%,17.56)
             (30\%,18.22)
          };
          \addplot  plot coordinates{

             (10\%,16.47)
             (20\%,17.99)
             (30\%,17.99)
          };

          \addplot  plot coordinates{

             (10\%,14.26)
             (20\%,17.47)
             (30\%,17.95)
          };
          \legend{\CDARP,\TCDARP,\CCDARP,\TCCDARP}
\end{axis}
\end{tikzpicture}
\caption{Savings from collaboration on group B instances}
\label{fig2}
\end{figure}

The savings in terms of overall objective function, using values of $\alpha_T=\alpha_C=\alpha$ ranging from 10\% to 30\%, are shown in Figures \ref{fig1} and \ref{fig2} for group A and B, respectively.
Reported values are averaged over the 28 instances of each group.
Savings obtained with \CDARP\  are reported as they represent an upper bound on the savings that can be obtained by the models.  On the other hand, the \TCCDARP\ model gives a lower bound  being the most constrained collaborative scenario.

The savings for group A are depicted in Figure \ref{fig1},
where the upper bound on the savings, achieved by the \CDARP, is 14.45\%.
Models \CCDARP, \TCDARP\ and \TCCDARP\ allow savings not less than 11.20\% with a loss not larger than 3.20\% with respect to \CDARP, when $\alpha=10\%$. When $\alpha$ increases to 30\%, savings improve up to 13.28\% with a loss with respect to \CDARP\ not larger than 1.17\% in the worst case.

The savings on group B instances are depicted in Figure \ref{fig2}. Here we see that with a slightly larger number of involved customers  we have more savings opportunities, and the upper bound on the savings 
becomes 18.51\%, while the \CCDARP\ and \TCDARP\ provide much closer performances: in the worst case the saving is just 2.39\% less than \CDARP\ for $\alpha=10\%$, and as little as $0.52\%$ less for $\alpha=30\%$.

\begin{table}[htbp!]\centering
\resizebox{\textwidth}{!}{
\begin{tabular}{|l|l|l|l|l|l|l|l|l|l|}
\hline
\multirow{3}{*}{Model}  & \multirow{3}{*}{\multirow{2}{*}{$\alpha$ (\%)}} & \multicolumn{4}{l|}{Group A}                                                                                           & \multicolumn{4}{l|}{Group B}                                                                                           \\ \cline{3-10}
                        &                                                                                    & \multirow{2}{*}{$\bar{U}$} & \multirow{2}{*}{$\widehat{U}$} & \multirow{2}{*}{$\bar{S}$} & \multirow{2}{*}{$\widehat{S} $} & \multirow{2}{*}{$\bar{U}$} & \multirow{2}{*}{$\widehat{U}$} & \multirow{2}{*}{$\bar{S}$} & \multirow{2}{*}{$\widehat{S} $}\\
                        &                                                                                    &                                          &                                         &                                            &                                  &                                          &                                         &                                             &                                      \\ \hline
\CDARP                  & -                                                                                  & 0.96                                     & 4                                       & 623.53                                     & 2262                       & 1.28                                     & 5                                       & 619.5                                       & 2014                                     \\ \hline
\TCDARP                & 10                                                                                 & 0.39                                     & 2                                       & 149.25                                     & 417                       & 0.82                                     & 3                                       & 200.54                                      & 421                                       \\ \cline{3-10}
                        & 20                                                                                 & 0.54                                     & 2                                       & 296.00                                     & 835                       & 1.00                                     & 3                                       & 373.14                                      & 1022                                   \\ \cline{3-10}
                        & 30                                                                                 & 0.64                                     & 2                                       & 363.96                                     & 921                       & 0.92                                     & 2                                       & 440.5                                       & 979                                       \\ \hline
\CCDARP                & 10                                                                                 & 0                                        & 0                                       & 215.96                                     & 600                       & 0.00                                     & 0                                       & 369.82                                      & 1218                                       \\ \cline{3-10}
                        & 20                                                                                 & 0                                        & 0                                       & 224.57                                     & 1076                      & 0.68                                     & 1                                       & 465.14                                      & 1639                                     \\ \cline{3-10}
                        & 30                                                                                 & 0.61                                     & 1                                       & 409.86                                     & 1478                      & 0.68                                     & 1                                       & 465.14                                      & 1639                                      \\ \hline
\TCCDARP               & 10                                                                                 & 0                                        & 0                                       & 167                                        & 417                       & 0.00                                     & 0                                       & 227.43                                      & 497                                       \\ \cline{3-10}
                        & 20                                                                                 & 0                                        & 0                                       & 202.29                                     & 645                       & 0.79                                     & 1                                       & 396.18                                      & 879                                      \\ \cline{3-10}
                        & 30                                                                                 & 0.61                                     & 1                                      & 364.5                                      & 921                        & 0.71                                     & 1                                       & 413.39                                      & 1160                                     \\ \hline
\end{tabular}}
\caption{Time and customers balance for group A and B instances}\label{invAB}
\end{table}


In Table \ref{invAB} statistics on the measures of balance, in terms of both time and customers, are reported.
It can be noted that the \CDARP\  produces a high level of exchange imbalance among companies.
Average and maximum time balance produced by the \CDARP\ for group A instances is 623.53 and 2262, respectively. For the \TCDARP\ the values are 363.96 and 921 in the worst case.
When we consider the customers balance, little can be said as, due to the limited size of the instances that can be solved to optimality, the number of customers is very low and at most 1 customer is available for exchange in the \CCDARP.
However, the \CDARP\ tends to exchange more customers between companies: 1 on average and up to 4  (50\% of companies' customers in group A).
A similar behaviour can be observed for group B.

Table \ref{invAB}, paired with the collaboration savings of Figures \ref{fig1} and \ref{fig2}, shows that it is possible to guarantee a time and/or customers balancing among companies without loosing  much in terms of collaboration savings.

\begin{figure}
\subfigure[$\alpha=10\%$]{\resizebox{0.45\textwidth}{!}{
\begin{tikzpicture}
\begin{axis}[xtick=data,legend style={at={(0.5,-0.15)},
      anchor=north,legend columns=-1},
    ,ybar,bar width=0.05,
    xlabel={$Days$},ylabel={\SAV\ (\%)}]

\addplot[fill=blue] coordinates {(1,35.94)(2,36.16)(3,32.38)(4,39.39)(5,32.28)(6,32.82)(7,32.51)
};
\addplot[fill=cyan] coordinates {(1,34.81)(2,35.89)(3,29.53)(4,36.99)(5,30.72)(6,32.11)(7,30.53)

};
\addplot[fill=lime] coordinates {(1,34.81)(2,35.89)(3,32.33)(4,39.34)(5,32.23)(6,32.77)(7,32.46)

};
\addplot[fill=magenta] coordinates {(1,34.81)(2,35.89)(3,29.53)(4,29.55)(5,30.72)(6,32.11)(7,30.53)

};

\addlegendentry{\CDARP}
\addlegendentry{\CCDARP}
\addlegendentry{\TCDARP}
\addlegendentry{\TCCDARP}

\end{axis}
\end{tikzpicture}}}
\quad\quad\subfigure[$\alpha=30\%$]{\resizebox{0.45\textwidth}{!}{
\begin{tikzpicture}
\begin{axis}[xtick=data,legend style={at={(0.5,-0.15)},
      anchor=north,legend columns=-1},
    ,ybar,bar width=0.05,
    xlabel={$Days$},ylabel={\SAV\ (\%)}]

\addplot[fill=blue] coordinates {(1,35.94)(2,36.16)(3,32.38)(4,39.39)(5,32.28)(6,32.82)(7,32.51)
};
\addplot[fill=cyan] coordinates {(1,35.89)(2,35.89)(3,32.33)(4,39.34)(5,30.72)(6,32.77)(7,32.46)

};
\addplot[fill=lime] coordinates {(1,35.89)(2,35.89)(3,32.33)(4,39.34)(5,32.23)(6,32.77)(7,32.46)

};
\addplot[fill=magenta] coordinates {(1,35.89)(2,35.89)(3,32.33)(4,39.34)(5,30.72)(6,32.11)(7,32.46)

};

\addlegendentry{\CDARP}
\addlegendentry{\CCDARP}
\addlegendentry{\TCDARP}
\addlegendentry{\TCCDARP}

\end{axis}
\end{tikzpicture}}}
\caption{Collaboration savings}\label{onegain}
\end{figure}
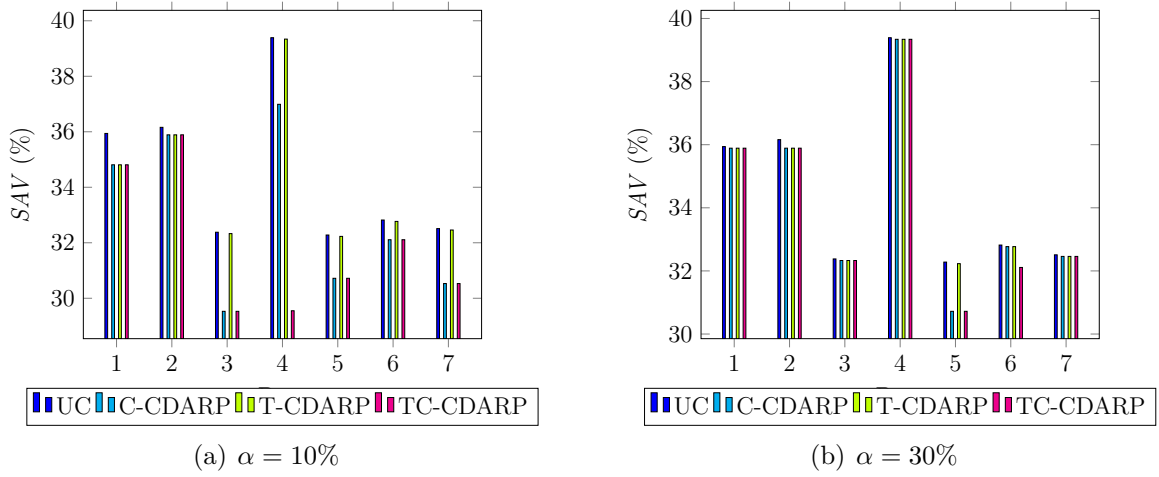

The same phenomenon can be observed in Figure \ref{onegain}, where savings are presented on the 7 instances of group B-R7, representing consecutive days,  and where the time and customers balance obtained in one day is respectively used to offset the balance of the subsequent day.
Namely, in day 1 we fix $S_m'=U_m'=0$, while for the following days the values of $S_m'$ and $U_m'$ are taken from the values of $S_m$ and $U_m$ of the day before.
The savings of model \TCDARP\ are quite close to the upper bound \CDARP, even considering the case of $\alpha=10\%$. While customers exchange on small size instances  in models \CCDARP\ and \TCCDARP\ is limited,  in most days the models obtain definitely good savings.
When $\alpha$ is set to 30\% all models perform almost as the bound \CDARP.

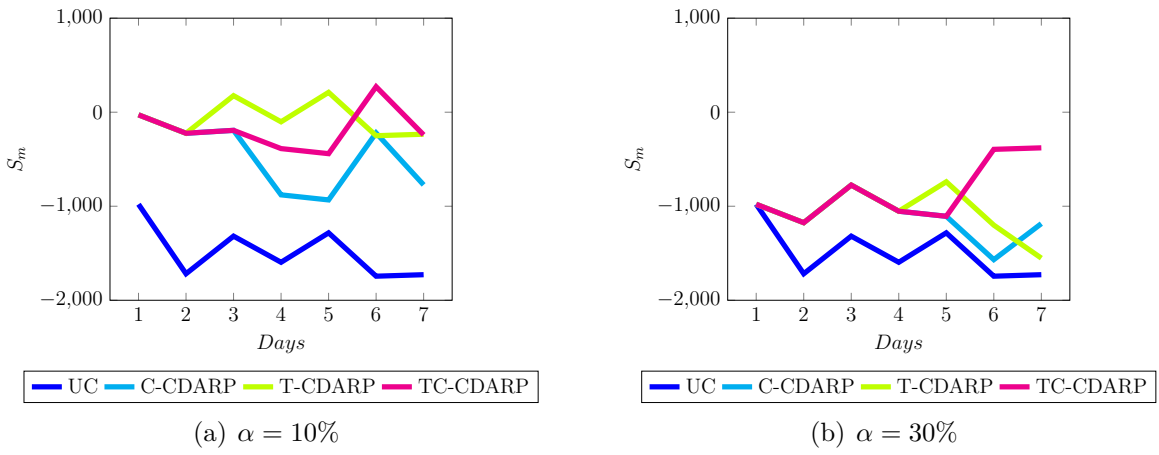
\begin{figure}
\subfigure[$\alpha=10\%$]{\resizebox{0.45\textwidth}{!}{
\begin{tikzpicture}
\begin{axis}[ymin=-2000, ymax = 1000
             ,xtick=data,legend style={at={(0.5,-0.25)},
      anchor=north,legend columns=-1},
    xlabel={$Days$},ylabel={$S_m$}]

\addplot[blue,line width=3pt] coordinates {(1,-979)(2,-1718)(3,-1318)(4,-1596)(5,-1283)(6,-1744)(7,-1728)

};
\addplot[cyan,line width=3pt] coordinates {(1,-28)(2,-224)(3,-193)(4,-879)(5,-933)(6,-221)(7,-773)

};
\addplot[lime,line width=3pt] coordinates {(1,-28)(2,-224)(3,176)(4,-102)(5,211)(6,-250)(7,-234)

};
\addplot[magenta,line width=3pt] coordinates {(1,-28)(2,-224)(3,-193)(4,-387)(5,-441)(6,271)(7,-241)

};

\addlegendentry{\CDARP}
\addlegendentry{\CCDARP}
\addlegendentry{\TCDARP}
\addlegendentry{\TCCDARP}

\end{axis}
\end{tikzpicture}}}
\quad\quad\subfigure[$\alpha=30\%$]{\resizebox{0.45\textwidth}{!}{
\begin{tikzpicture}
\begin{axis}[ymin=-2000, ymax = 1000
             ,xtick=data,legend style={at={(0.5,-0.25)},
      anchor=north,legend columns=-1},
    xlabel={$Days$},ylabel={$S_m$}]

\addplot[blue,line width=3pt] coordinates {(1,-979)(2,-1718)(3,-1318)(4,-1596)(5,-1283)(6,-1744)(7,-1728)

};
\addplot[cyan,line width=3pt] coordinates {(1,-979)(2,-1175)(3,-775)(4,-1053)(5,-1107)(6,-1568)(7,-1185)

};
\addplot[lime,line width=3pt] coordinates {(1,-979)(2,-1175)(3,-775)(4,-1053)(5,-740)(6,-1201)(7,-1552)

};
\addplot[magenta,line width=3pt] coordinates {(1,-979)(2,-1175)(3,-775)(4,-1053)(5,-1107)(6,-395)(7,-379)

};

\addlegendentry{\CDARP}
\addlegendentry{\CCDARP}
\addlegendentry{\TCDARP}
\addlegendentry{\TCCDARP}

\end{axis}
\end{tikzpicture}}}
\caption{Time balance of the company in 7 days}\label{onetime}
\end{figure}

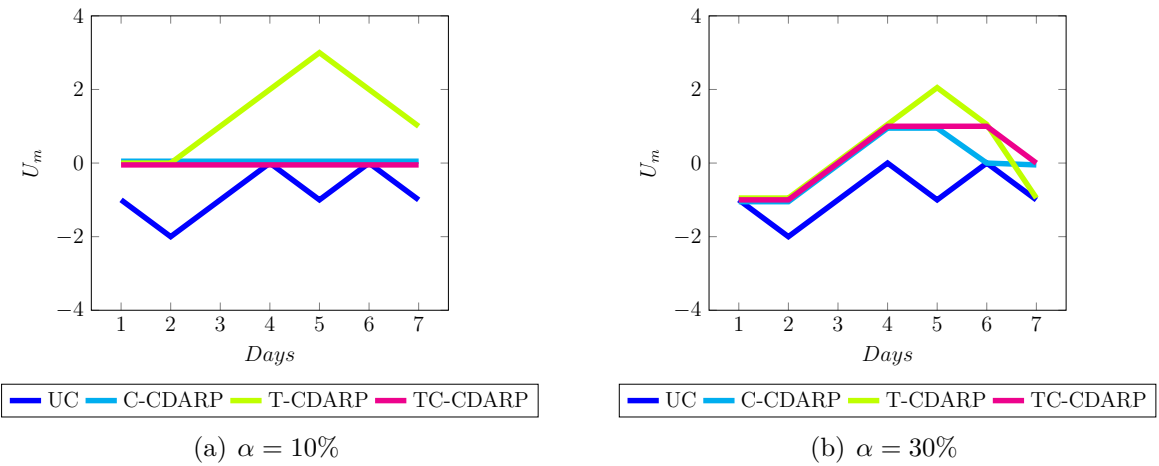
\begin{figure}
\subfigure[$\alpha=10\%$]{\resizebox{0.45\textwidth}{!}{
\begin{tikzpicture}
\begin{axis}[ymin=-4, ymax = 4
             ,xtick=data,legend style={at={(0.5,-0.25)},
      anchor=north,legend columns=-1},
    xlabel={$Days$},ylabel={$U_m$}]

\addplot[blue,line width=3pt] coordinates {(1,-1)(2,-2)(3,-1)(4,0)(5,-1)(6,0)(7,-1)
};
\addplot[cyan,line width=3pt] coordinates {(1,0.05)(2,0.05)(3,0.05)(4,0.05)(5,0.05)(6,0.05)(7,0.05)

};
\addplot[lime,line width=3pt] coordinates {(1,0)(2,0)(3,1)(4,2)(5,3)(6,2)(7,1)

};
\addplot[magenta,line width=3pt] coordinates {(1,-0.05)(2,-0.05)(3,-0.05)(4,-0.05)(5,-0.05)(6,-0.05)(7,-0.05)

};

\addlegendentry{\CDARP}
\addlegendentry{\CCDARP}
\addlegendentry{\TCDARP}
\addlegendentry{\TCCDARP}

\end{axis}
\end{tikzpicture}}}
\quad\quad\subfigure[$\alpha=30\%$]{\resizebox{0.45\textwidth}{!}{
\begin{tikzpicture}
\begin{axis}[ymin=-4, ymax = 4
             ,xtick=data,legend style={at={(0.5,-0.25)},
      anchor=north,legend columns=-1},
    xlabel={$Days$},ylabel={$U_m$}]

\addplot[blue,line width=3pt] coordinates {(1,-1)(2,-2)(3,-1)(4,0)(5,-1)(6,0)(7,-1)
};
\addplot[cyan,line width=3pt] coordinates {(1,-1.05)(2,-1.05)(3,-0.05)(4,0.95)(5,0.95)(6,0)(7,-0.05)

};
\addplot[lime,line width=3pt] coordinates {(1,-0.95)(2,-0.95)(3,0.05)(4,1.05)(5,2.05)(6,1.05)(7,-0.95)

};
\addplot[magenta,line width=3pt] coordinates {(1,-1)(2,-1)(3,0)(4,1)(5,1)(6,1)(7,0)

};

\addlegendentry{\CDARP}
\addlegendentry{\CCDARP}
\addlegendentry{\TCDARP}
\addlegendentry{\TCCDARP}

\end{axis}
\end{tikzpicture}}}
\caption{Customer balance of the company in 7 days}\label{onecust}
\end{figure}

In Figure \ref{onetime} we present the time  balance of one company. We see that when $\alpha=10\%$, the \TCDARP\ and the \TCCDARP\  are able to maintain the time balance low through the seven days as the \CCDARP, to a lesser extent, also does. On the contrary, the \CDARP\ produces remarkable imbalances. The desired effect of using constrained models is clearly reduced when using $\alpha=30\%$ as the bounds on time and customers balance are less tight.
In Figure \ref{onecust}, the customers balance is depicted. As expected, the \CCDARP\ and the \TCCDARP\  are able to maintain the customer balance low over days with both $\alpha=10\%$ and $\alpha=30\%$. The impact of the \TCDARP\  on the customer balance is low, similarly to the impact of the \CDARP.


\subsection{Assessment  of the ALNS performance against optimal solutions}\label{oneinstance}

 Here we  report a comparison of the performance of the ALNS heuristic with respect to the optimal solutions obtained on group A and B instances.
In the computational experiments,  the optimal solution of the non-collaborative setting was used as initial solution for the heuristic.

Table \ref{time} shows the computational time needed by CPLEX to solve the \CDARP, the \TCDARP, the \CCDARP\ and the \TCCDARP\  along with the ALNS computational time.
We see that, while CPLEX is faster than the ALNS heuristic on group A instances,  as the size increases from 8 to 10 customers the computational times for CPLEX widely increase.  When the size is further increased to more than two companies or/and more than 10 requests, CPLEX is not able to provide a solution within 7200 seconds.
As expected, the ALNS computational times are few seconds for both groups A and B.

\begin{table}[htbp!]\centering
\begin{tabular}{|l|l|ll|ll|}
\hline
\multirow{2}{*}{Model} & \multirow{2}{*}{$\alpha$ (\%)} & \multicolumn{2}{c|}{Group A (sec.)}            & \multicolumn{2}{c|}{Group B (sec.)}           \\
                        &                            & \multicolumn{1}{l|}{CPLEX} & ALNS & \multicolumn{1}{l|}{CPLEX} & ALNS \\ \hline
\CDARP     & -                          & 39.6                            & 12.2      & 1692.7                          & 14.6      \\ \hline
\TCDARP                   & 10                        & 2.8                             & 12.9      & 90.7                            & 13.6      \\
                        & 20                        & 3.7                             & 14.4      & 141.3                          & 13.9      \\
                        & 30                        & 4.4                             & 13.0      & 206.6                          & 14.2      \\ \hline
\CCDARP                & 10                        & 2.8                             & 12.7      & 61.2                           & 12.8      \\
                        &20                        & 2.8                             & 12.5      & 119.1                           & 13.7      \\
                        & 30                        & 4.6                             & 11.0      & 119.1                           & 14.1      \\ \hline
\TCCDARP       & 10                        & 2.2                             & 13.7      & 48.1                            & 13.6      \\
                        & 20                        & 2.5                             & 12.7      & 106.6                           & 13.5      \\
                        & 30                        & 4.0                             & 12.3      & 124.7                           & 13.5      \\ \hline
\end{tabular}
\caption{Computational times on group A and B instances}
\label{time}
\end{table}

In Table \ref{gap}, the quality of the solutions produced by the ALNS heuristic  is assessed by comparison to the optimal ones produced by CPLEX on instances of groups A and B.
Results for each model are grouped over different values of parameter $\alpha$. The table shows the ALNS average and maximum optimality gap, the \% of instances solved to optimality, and the percentage of instances not solved to optimality in error classes with step 1\%. The ALNS solves to optimality 85\% of instances in groups A and B and 94.5\% of the instances are less than 1\% away from the optimum.
%


%

\begin{table}[htbp!]\centering
\resizebox{\textwidth}{!}{\begin{tabular}{|l|l|l|l|l|l|l|l|l|l|l|l|l|l|}
\hline
\multirow{2}{*}{Model} & \multicolumn{7}{c|}{\GAP\ on Group A}                     & \multicolumn{6}{c|}{\GAP\ on Group B}                    \\
                        &                                                                                    Average & Max  & opt.\% & 0-1\% & 1-2\% & 2-3\% & 3-4\% & Average & Max  & opt.\% & 0-1\% & 1-2\% & 2-3\%  \\ \hline

                        \CDARP                                                                                                   & 0.05    & 1.49 & 96.4 & 0.0   & 3.6   & 0.0   & 0.0  &
                        0.04& 0.81 & 89.3 & 10.7  & 0.0   & 0.0     \\ \hline
\TCDARP                                                                                                & 0.13    & 2.41 & 85.6 & 7.2   & 4.8   & 2.4   & 0.0  &

0.22    & 2.57 & 77.4 & 14.3  & 5.9   & 2.4    \\ \hline
\CCDARP                                                                                                & 0.26    & 3.6 & 89.2 & 2.4   & 0.0   & 3.6   & 4.8  &

0.16    & 2.83 & 78.6 & 17.9  & 2.4   & 1.2   \\ \hline
\TCCDARP                                                                                               & 0.16    & 3.32 & 85.7 & 9.5   & 0.0   & 2.4   & 2.4  &

0.23    & 2.69 & 76.2 & 15.5  & 4.8  & 3.6   \\ \hline
\end{tabular}}
\caption{Performance of the ALNS heuristic}\label{gap}
\end{table}

\subsection{Savings on large scale instances}\label{allinstances}

The ALNS heuristic has been run to solve models \NCDARP, \CDARP, \TCDARP, \CCDARP\ and \TCCDARP\ on instances of groups C e D, with $\alpha$ ranging from 10\% to 30\%. Results are averaged over the 56 instances of groups C and  D.

\begin{figure}[htbp!]
\centering
\begin{tikzpicture}
\begin{axis} [width=8cm,
              symbolic x coords={10\%,20\%,30\%},xtick=data,
legend style={at={(0.5,1.15)},
	anchor=north,legend columns=-1},ymin=0,
    xlabel={$\alpha$},ylabel={\SAV\ (\%)}]
\addplot  plot coordinates{

             (10\%,32.12)
             (20\%,32.12)
             (30\%,32.12)
          };
            \addplot  plot coordinates{

             (10\%,29.20)
             (20\%,31.02)
             (30\%,31.40)
          };
          \addplot  plot coordinates{

             (10\%,27.62)
             (20\%,30.20)
             (30\%,30.71)
          };

          \addplot  plot coordinates{

             (10\%,27.18)
             (20\%,29.99)
             (30\%,30.50)
          };
          \legend{\CDARP,\TCDARP,\CCDARP,\TCCDARP}
\end{axis}
\end{tikzpicture}
\caption{Collaboration saving on group C instances}\label{collgainHC}
\end{figure}
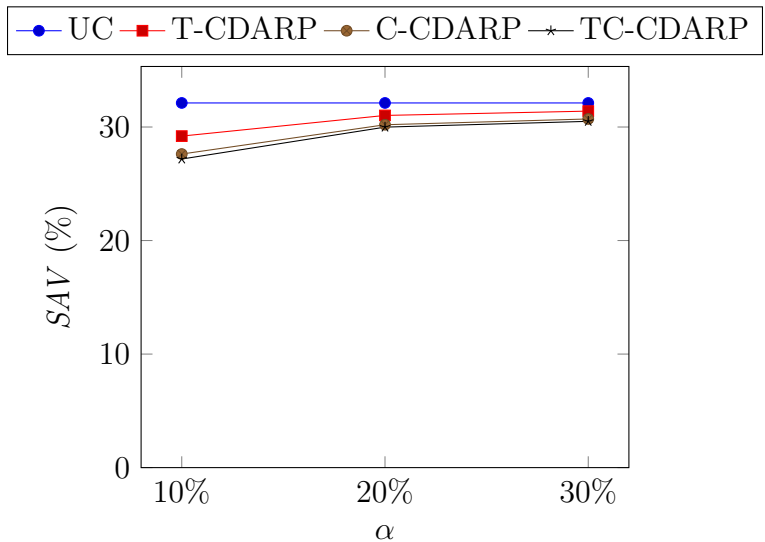

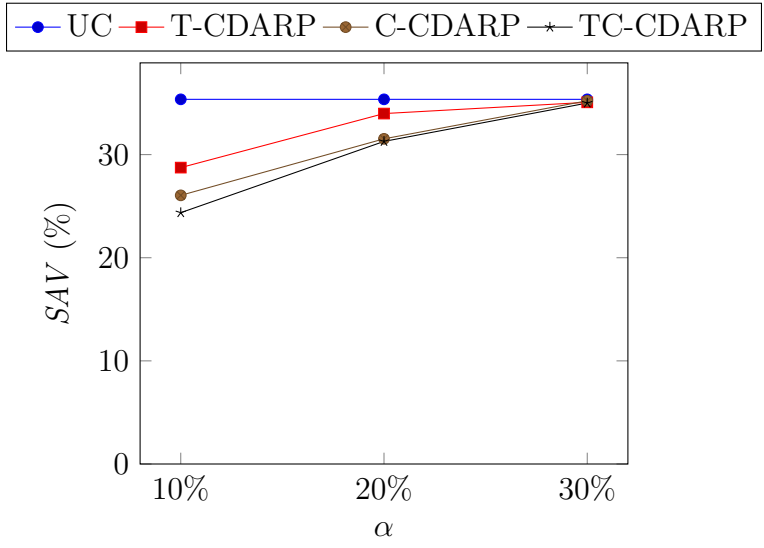
\begin{figure}[htbp!]
\centering
\begin{tikzpicture}
\begin{axis} [width=8cm,
              symbolic x coords={10\%,20\%,30\%},xtick=data,
legend style={at={(0.5,1.15)},
	anchor=north,legend columns=-1},ymin=0,
    xlabel={$\alpha$},ylabel={\SAV\ (\%)}]
\addplot  plot coordinates{

             (10\%,35.36)
             (20\%,35.36)
             (30\%,35.36)
          };
           \addplot  plot coordinates{

             (10\%,28.74)
             (20\%,33.98)
             (30\%,35.07)
          };
          \addplot  plot coordinates{

             (10\%,26.05)
             (20\%,31.53)
             (30\%,35.22)
          };

          \addplot  plot coordinates{

             (10\%,24.37)
             (20\%,31.28)
             (30\%,35.00)
          };
          \legend{\CDARP,\TCDARP,\CCDARP,\TCCDARP}
\end{axis}
\end{tikzpicture}
\caption{Collaboration saving on group D instances}\label{collgainHD}
\end{figure}

The savings, using values of $\alpha_T = \alpha_C = \alpha$
ranging from 10\% to 30\%, are shown in Figures \ref{collgainHC} and \ref{collgainHD} for group C and D, respectively. Savings obtained with
\CDARP\ are reported as they represent an upper bound on the savings that can be obtained by the collaborative models.
The savings for group C instances  are depicted in Figure \ref{collgainHC}, where the upper bound on the savings, obtained with the \CDARP, is 32.12\%.
All the models show an excellent performance as, even with the most constrained \TCCDARP\ model, the difference in terms of savings with respect to the upper bound \CDARP is  small, less than 5\% for $\alpha=10\%$ and less than 1.63\% for $\alpha=30\%$.
In Figure \ref{collgainHD}, we observe that increasing the number of customers the saving of the \CDARP\  with respect to the non-collaborative setting increases to 35.36\%.
On the other side we have a loss in savings up to 11\% in the most constrained case for $\alpha=30\%$. The losses decrease with $\alpha$ and become negligible for $\alpha=30\%$.

In Table \ref{invCD}, statistics on the balance of both times and customers with respect to the non-collaborative case are shown.
Also here, a high level of imbalance can be observed when an unconstrained setting is adopted.
The imbalance in customers and time is strongly mitigated by applying the \CCDARP\ and the \TCDARP, respectively. The two issues are simultaneously fixed by the \TCCDARP\ model.
Table \ref{invCD}, paired with the savings shown in Figures \ref{collgainHC} and \ref{collgainHD}, shows that it is possible to guarantee a time and/or customers balance among companies without loosing  much in terms of collaboration saving also on group C and D instances.

The computational time needed to run the ALNS heuristic on these instances is very small, as shown in Table \ref{timeH}.
We get slightly higher times for $\alpha=10\%$, that is when the possibilities to exchange customers are more constrained.


\begin{table}[htbp!]\centering
\resizebox{\textwidth}{!}{
\begin{tabular}{|l|l|l|l|l|l|l|l|l|l|}
\hline
\multirow{3}{*}{Model} & \multirow{3}{*}{\multirow{2}{*}{$\alpha$ (\%)}} & \multicolumn{4}{l|}{Group C}                                                                                           & \multicolumn{4}{l|}{Group D}                                                                                           \\ \cline{3-10}
                        &                                                                                    & \multirow{2}{*}{$\bar{U}$} & \multirow{2}{*}{$\widehat{U}$} & \multirow{2}{*}{$\bar{S}$} & \multirow{2}{*}{$\widehat{S} $} &
                        \multirow{2}{*}{$\bar{U}$} & \multirow{2}{*}{$\widehat{U}$} & \multirow{2}{*}{$\bar{S}$} & \multirow{2}{*}{$\widehat{S} $} \\
                        &                                                                                    &                                          &                                         &                                            &                                       &                                           &                                         &                                             &                               \\ \hline
\CDARP                  & -                                                                                  & 1.61                                     & 5                                       & 778.84                                     & 3181                                  & 1.59                                      & 8                                       & 812.32                                      & 3219                                    \\ \hline
\TCDARP                & 10                                                                                 & 1.23                                     & 5                                       & 527.36                                     & 1442                                     & 1.12                                     & 5                                       & 369.05                                      & 1112                                     \\ \cline{3-10}
                        & 20                                                                                 & 1.67                                     & 5                                       & 595.46                                     & 1690                                      & 1.38                                     & 6                                       & 668.11                                      & 2268                                    \\ \cline{3-10}
                        & 30                                                                                 & 1.45                                     & 5                                       & 692.57                                     & 1920                                      & 1.59                                     & 6                                       & 775.01                                      & 3031                                    \\ \hline
\CCDARP                & 10                                                                                 & 0.62                                     & 1                                       & 484.32                                     & 2159                                  & 0.58                                     & 1                                       & 423.77                                     & 1765                                    \\ \cline{3-10}
                        & 20                                                                                 & 1.02                                     & 2                                       & 601.57                                     & 2435                                    & 1.13                                     & 2                                       & 599.82                                     & 2455                                     \\ \cline{3-10}
                        & 30                                                                                 & 1.23                                     & 3                                       & 632.82                                     & 2146                                     & 1.40                                     & 3                                       & 715.71                                      & 2402                                     \\ \hline
\TCCDARP               & 10                                                                                 & 0.59                                     & 1                                       & 401.09                                     & 1208                                   & 0.55                                     & 1                                       & 297.18                                     & 1037                                     \\ \cline{3-10}
                        & 20                                                                                 & 1.11                                     & 2                                       & 595.46                                     & 1690                                      & 1.05                                     & 2                                       & 588.03                                      & 2135                                     \\ \cline{3-10}
                        & 30                                                                                 & 1.36                                     & 3                                       & 692.57                                     & 1920                                   & 1.42                                     & 3                                       & 696.17                                      & 2552                                     \\ \hline
\end{tabular}}
\caption{Time and customer  deviation for group C and D instances}\label{invCD}
\end{table}

\begin{table}[htbp!]
\centering
\begin{tabular}{|l|l|l|l|}
                    \hline   & \multicolumn{3}{|l|}{Computational time (sec)} \\\hline
Model                & $\alpha$ (\%)       & group C       & group D      \\ \hline
\CDARP
                      & -            & 7.86        & 23.76        \\ \hline
\TCDARP                  & 10            & 8.26        & 26.45        \\
                      & 20            & 6.65        & 23.58        \\
                      & 30            & 7.56        & 23.33        \\ \hline
\CCDARP              & 10            & 8.14        & 26.61        \\
                      & 20            & 6.69        & 23.32        \\
                      & 30            & 6.83        & 23.17        \\ \hline
\TCCDARP   & 10            & 14.24        & 31.02        \\
                      & 20            & 7.19        & 24.14        \\
                      & 30            & 7.82        & 26.07     \\ \hline
\end{tabular}
\caption{Computational time for group C and  D instances}\label{timeH}
\end{table}

\subsection{Destroy/repair operator effectiveness}\label{operatoreffectiveness}

In Table \ref{desCDHrank} we show the aggregated impact of destroy operators on all groups of instances.
The results are  aggregated also with respect to  $\alpha$ as we noticed no relevant insights on detailed results.

Each row gives the ranking of the operators, in terms of number of times in which the different destroy operators have been successful in finding an improving solution for a given model. %
For instances of groups A, B and C, the Closeness removal and the Interchangeability removal are the most successful operators for all models while Worst and Related are the least successful ones. In group D instances, the Random and the Related operators rank first and second, respectively. It is worth observing that, overall, the Closeness and the Interchangeability are the most successful operators.
In Table \ref{repCDHrank} the same statistic is reported for the repair operators.
The 2-Regret operator is the most successful in groups A and B instances whereas the operator 3-Regret ranks second for group A instances and the Closeness insertion ranks second for group B instances.
The operator 4-Regret is the least successful on both groups.
Conversely to the results obtained on groups A and B instances, operator 4-Regret ranks second in group C and D instances. 2-Regret operator remains the most successful in group C instances while the operator Best insertion ranks first in group D instances. 3-Regret operator seems to perform slightly worse than the other operators for both group C and D instances. Overall, the operators 2-Regret and Closeness insertion rank, respectively, first and second most successful ones.

We have run all the experiments on groups C and D instances also removing all the new destroy and repair operators, namely the Closeness removal and  the Interchangeability removal operators and the Closeness insertion operator, to test their level of relevance. In fact, using those operators reduces the objective function on average of the 0.2\% for group C instances and of  1.06\% for group D instances, which assesses their positive impact.

\begin{table}[htbp!]\centering
\resizebox{\textwidth}{!}{
\begin{tabular}{|l|l|l|l|l|l|l|l|}
\hline
      Group          & Model & Random & Worst & Related & Proximity & Closeness & Interchangeability \\ \hline
      \multirow{3}{*}{A}
                 & \TCDARP       & IV& VI& V& III& I& II

 \\
              & \CCDARP      & IV& VI& V& III& I& II

 \\
   & \TCCDARP       & IV& VI& V& III& I& II
\\\cline{2-8}
   & All models       & IV& VI& V& III& I& II

 \\ \hline
 \multirow{4}{*}{B}& \TCDARP       & III& VI& IV& V& I& II

 \\
              & \CCDARP      & II& VI& V& IV& I& III

 \\
   & \TCCDARP       & V& VI& IV& II& I& III
\\\cline{2-8}
   & All models       & III& VI& V& IV& I&II
 \\ \hline
 \multirow{4}{*}{C}& \TCDARP       & IV&VI&II&V&I&III
 \\
              & \CCDARP      & IV&VI&II&V&I&III
 \\
   & \TCCDARP       & III&V&VI&IV&II&I
   \\\cline{2-8}
   & All models       & IV & VI& III& V& I& II
 \\ \hline
       \multirow{4}{*}{D}& \TCDARP       &II&III&I&VI&V&IV
    \\
        & \CCDARP       & IV&V&I&VI&III&II
    \\
 & \TCCDARP       & V&VI&I&IV&II&III
 \\\cline{2-8}
   & All models       & I & V & II & VI & IV & III
 \\\hline
   All & All models       & III&VI&V&IV&\textbf{I}&\textbf{II}
   \\\hline
\end{tabular}}
\caption{Ranking of destroy operators }\label{desCDHrank}
\end{table}

\begin{table}[htbp!]\centering
\resizebox{\textwidth}{!}{
\begin{tabular}{|l|l|l|l|l|l|l|}
\hline
Group                & Model & Best insertion  & 2-Regret & 3-Regret & 4-Regret & Closeness\\ \hline
\multirow{4}{*}{A}& \TCDARP       & IV& I& III& V& II

 \\
              & \CCDARP      & III& I& II& IV& V

 \\
   & \TCCDARP       & III& I& IV& V& II
\\\cline{2-7}
   & All models       & IV& I& II& V& III
 \\ \hline
 \multirow{4}{*}{B}& \TCDARP       & III& II& V& IV& I

 \\
              & \CCDARP      & II& I& III& V& IV

 \\
   & \TCCDARP       & IV& I& III& V& II
\\\cline{2-7}
   & All models       & III & II & IV& V & I
 \\ \hline
\multirow{4}{*}{C}& \TCDARP       & III&II&V&I&IV
 \\
              & \CCDARP      & IV&I&V&III&II
 \\
   & \TCCDARP       & III&II&V&I&IV
   \\\cline{2-7}
   & All models       & III& I& V& II& IV
 \\ \hline
       \multirow{4}{*}{D} & \TCDARP  &     I&III&IV&II&V
    \\
        & \CCDARP       & III&I&V&II&IV
    \\
 & \TCCDARP       & I&III&V&II&IV
 \\ \cline{2-7}
   & All models       & I& III& V& II& IV
 \\\hline
 All& All models       & III&\textbf{I}&IV&V&\textbf{II}
 \\\hline
\end{tabular}}
\caption{Ranking of repair operators }\label{repCDHrank}
\end{table}

\section{Conclusions}\label{conclusions}

In this paper we have
studied the problem of optimizing the operations of companies, offering demand-responsive transportation services, involved in horizontal collaboration. We have proposed a collaboration scheme where each  company is allowed to serve customers of other companies, if beneficial. The exchange  is, however, constrained to make the collaboration scheme acceptable to all companies, in particular to those of small or relatively small size. This is achieved by imposing, in an optimization model, that the workload and/or the  number of customers served by each company does not change much with respect to the case where no collaboration scheme is put in place. An ALNS heuristic has been presented for the solution of the resulting mixed integer models that optimize the total traveling cost.

The computational experiments have shown that the ALNS heuristic has a good performance, with an average error of 0.3\%. More importantly, the optimization models allow a saving, with respect to a non-collaborative setting, of up to 30\%, very close to the upper bound on the possible savings that is obtained through the solution of an unconstrained collaboration model.

Several research directions remain to be explored. Horizontal collaboration is a hot topic nowadays, especially in transportation. It would be interesting to extend the proposed collaboration setting to other kinds of companies, in transportation but also in other domains.  Also, dynamic and real-time versions of the proposed approach would be worth of investigation.

\newpage
\bibliographystyle{natbib}
\bibliography{darpBib}
\end{document}